\documentclass[11pt,reqno]{amsart}

\usepackage{amssymb,amsmath,amsthm}
\usepackage{hyperref}
\usepackage{geometry}
\newtheorem{theorem}{Theorem}[section]
\newtheorem{proposition}[theorem]{Proposition}
\newtheorem{lemma}[theorem]{Lemma}
\newtheorem{corollary}[theorem]{Corollary}
\theoremstyle{definition}
\newtheorem{definition}[theorem]{Definition}
\newtheorem{remark}[theorem]{Remark}
\newtheorem{example}[theorem]{Example}

\numberwithin{equation}{section}

\begin{document}

\title[Complement Cancellation and the Spectral Characterization of Three-Point Recurrence]{Complement Cancellation and the Spectral Characterization of Three-Point Recurrence}

\author{Dibyendu De}
\address{Department of Mathematics, University of Kalyani, Kalyani 741235, West Bengal, India}
\email{dibyendude@klyuniv.ac.in}


\begin{abstract}
In this paper, we investigate both two-sided and one-sided \(2\)-nice recurrence within the framework of weakly mixing measure-preserving systems. As our principal result, we establish an exact spectral characterization for the failure of two-sided \(2\)-nice recurrence: namely, a subset \(R\subseteq\mathbb N\) fails to possess this property if and only if there exist a continuous symmetric probability measure \(\sigma\) on \(\mathbb T\) and a constant \(\varepsilon>0\) satisfying
\[
2\widehat{\sigma}(n)+\widehat{\sigma}(2n)\leq-\varepsilon
\qquad(n\in R).
\]
Our proof of the necessity relies upon a complement-symmetrized correlation identity, whereas the converse is achieved through a construction involving a stationary Gaussian system, a bounded odd observable, Hermite expansions, and an independent Bernoulli extension. Furthermore, we extend this characterization to three-point polynomial configurations as well as to actions of countable abelian groups. In addition, we derive an unconditional analogue in which arbitrary symmetric probability measures take the place of continuous ones, thereby demonstrating a strict separation between the ergodic and weakly mixing classes of recurrence.

Subsequently, we turn our attention to the distinction between two-sided and ordinary one-sided \(2\)-nice recurrence. It is shown that a suitable negative spectral condition compels the failure of one-sided recurrence. More remarkably, we construct an explicit weakly mixing system wherein the quadratic spectral expression remains uniformly positive, yet the presence of a negative cubic correlation still induces a one-sided recurrence deficit. Consequently, it follows that two-point spectral data alone are insufficient to govern one-sided \(2\)-nice recurrence.
\end{abstract}

\maketitle

\section{Introduction}

The study of recurrence in measure-preserving dynamical systems serves as a natural confluence of ergodic theory, harmonic analysis, and combinatorial number theory. Ever since Furstenberg's celebrated ergodic proof of Szemer?(C)di's theorem \cite{Fur77,Fur81}, a central theme in this area has been the intimate correspondence between measure-theoretic recurrence and combinatorial largeness. For single recurrence, spectral theory furnishes an especially transparent description: the two-point correlations of a function are given directly by the Fourier--Stieltjes coefficients of its spectral measure. Multiple recurrence, however, is considerably more delicate. Its correlation expressions involve several iterates of the transformation, and their analysis frequently demands sophisticated machinery such as characteristic nilfactors and PET induction \cite{HK05,Ber87}.

Let us recall the relevant notion of nice recurrence. Let \((X,\mathcal B,\mu,T)\) be an invertible probability-preserving system. For an integer \(k\geq 1\), a subset \(R\subseteq\mathbb N\) is termed a \emph{set of \(k\)-nice recurrence} for a specified class of systems if, for every system belonging to that class, every measurable set \(A\) with \(\mu(A)>0\), and every \(\varepsilon>0\), there exists \(n\in R\) such that
\[
 \mu\left(\bigcap_{j=0}^{k}T^{-jn}A\right)
 \geq \mu(A)^{k+1}-\varepsilon.
\]
The case \(k=1\) corresponds to the classical notion of nice recurrence. Recently, Chapman \cite{Cha26} established the partition regularity of sets of nice recurrence, while related developments concerning countable groups and inverse correspondence principles have appeared in \cite{FR26,FS25,RM25,FT26}. These investigations naturally motivate the question pursued in the present work: what remains of the spectral and combinatorial structure when one passes from single recurrence to three-point and higher-point configurations?

The principal object of study in this paper is a complement-symmetrized variant of \(2\)-nice recurrence. For a measurable set \(A\), let us write
\[
 b_n(A)=\mu(A\cap T^{-n}A\cap T^{-2n}A).
\]
We seek an integer \(n\in R\) satisfying
\[
 b_n(A)+b_n(A^c)
 \geq \mu(A)^3+\mu(A^c)^3-\varepsilon.
\]
This formulation is distinct from ordinary one-sided \(2\)-nice recurrence. Its decisive advantage lies in an exact algebraic cancellation: if we decompose \(1_A=\alpha+f\), where \(\alpha=\mu(A)\) and \(\int f\,d\mu=0\), then the centered expansion of \(b_n(A)\) involves a cubic term, whereas the corresponding expansion for \(A^c\) carries its exact negative. Upon adding these two expressions, the cubic terms cancel out, leaving solely the quadratic correlations \(r_n=\langle f,T^nf\rangle\) and \(r_{2n}\).

Our first main theorem elevates this elementary observation into a complete spectral characterization. Specifically, a subset \(R\subseteq\mathbb N\) fails two-sided \(2\)-nice recurrence for weakly mixing systems if and only if there exist a continuous symmetric probability measure \(\sigma\) on \(\mathbb T=\mathbb R/\mathbb Z\) and a constant \(\varepsilon_2>0\) such that
\[
 2\widehat{\sigma}(n)+\widehat{\sigma}(2n)
 \leq-\varepsilon_2\qquad(n\in R).
\]
The forward implication is obtained immediately by passing to the normalized spectral measure of \(1_A-\mu(A)\). For the converse, we realize \(\sigma\) as the spectral measure of a stationary Gaussian process. A bounded odd function of the coordinate Gaussian variable possesses a vanishing centered three-point correlation, while its leading Hermite component governs its two-point correlations. An independent Bernoulli extension then converts the resulting \([0,1]\)-valued observable into the indicator function of a measurable set without altering the pertinent multiple correlations. This mechanism also clarifies the precise role played by weak mixing in the characterization: it corresponds exactly to the continuity of the spectral measure. Furthermore, we demonstrate that when all invertible probability-preserving systems are admitted, the same characterization persists without any continuity requirement on the symmetric spectral measure. An argument utilizing positive upper Banach density subsequently yields an explicit separation between the weakly mixing and ergodic recurrence classes.

The very same algebraic mechanism extends seamlessly to general configurations of the form \(\{0,p_1(n),p_2(n)\}\), where the associated quadratic spectral expression takes the form
\[
 \widehat{\sigma}(p_1(n))+
 \widehat{\sigma}(p_2(n))+
 \widehat{\sigma}((p_2-p_1)(n)).
\]
We formulate this extension both for integer-valued polynomial configurations and, more broadly, for suitable homomorphism configurations in countable abelian groups.

Nonetheless, the exact two-sided characterization does not resolve the one-sided problem. Indeed, the identity
\[
 b_n(A)-\mu(A)^3
 =\mu(A)(2r_n+r_{2n})+c_n,
 \qquad
 c_n=\int f(T^nf)(T^{2n}f)\,d\mu,
\]
reveals that the cubic correlation remains operative when considering a single set. We examine this distinction by employing Gaussian observables incorporating both odd and even Hermite components. The resulting construction is carefully designed to exhibit a strictly positive quadratic spectral contribution alongside a negative cubic correlation. The complete proof and the choice of parameters are presented in Section~\ref{sec:onesided}.

The remainder of the paper is organized as follows. In Section~2, we assemble the necessary measure-theoretic, spectral, and Gaussian preliminaries. Section~3 is devoted to introducing the recurrence notions, establishing the exact spectral characterization, and deriving its polynomial and group-theoretic extensions. Finally, in Section~4, we investigate the one-sided problem and construct the cubic-driven counterexample outlined above.

\bigskip

\noindent\textbf{AI use acknowledgement: }Chat GPT was used during early brainstorming and for language editing. All mathematical arguments were independently developed, checked, and verified by the author, who assumes full responsibility for the contents.
\section{Preliminaries}
Throughout this paper, \(\mathbb{N}=\{1,2,\dots\}\) denotes the set of positive integers and \(\mathbb{T}=\mathbb{R}/\mathbb{Z}\) denotes the circle group. We write \(e(t)=e^{2\pi it}\) for \(t\in\mathbb{T}\). Unless explicitly stated otherwise, all measure-preserving systems considered herein are assumed to be invertible probability-preserving systems.

The objective of this section is to bring together the dynamical and harmonic-analytic tools required for our investigation of two-sided \(2\)-nice recurrence. A guiding principle in this discipline is that recurrence properties of subsets of the integers can frequently be illuminated by passing between combinatorial configurations and correlations in measure-preserving systems. Bergelson's density theorem generalizing Schur's theorem \cite{Ber86Schur} affords an instructive illustration: an ergodic-theoretic recurrence argument is translated into a combinatorial density result via finite averaging and intersection estimates.

In the present context, we pursue a more specialized question. For a measurable set \(A\), the three-term correlation involves both the two-point correlations of its centered indicator function and a genuinely three-point correlation. Our crucial observation is that passing simultaneously to \(A\) and its complement \(A^c\) cancels the latter term, leaving behind a quadratic expression involving correlations at times \(n\) and \(2n\). Spectral measures thus emerge as the natural language for analyzing this phenomenon.

Accordingly, we first review the fundamental notions of measure-preserving and weakly mixing systems, together with the spectral representation of correlations. We then introduce stationary Gaussian systems and Hermite expansions, which furnish the converse mechanism: any prescribed continuous symmetric spectral measure can be realized dynamically, while a suitably chosen bounded observable allows its first Hermite component to dominate the relevant quadratic correlations. An independent Bernoulli extension will subsequently convert this functional witness into an explicit measurable set. Thus, the material gathered here is not merely introductory; it provides the essential framework enabling the passage from an abstract Fourier inequality to a concrete failure of recurrence.

\subsection{Measure-preserving systems and recurrence}
A \emph{probability-preserving system} is a quadruple \((X,\mathcal{B},\mu,T)\), where \((X,\mathcal{B},\mu)\) is a standard probability space and \(T\colon X\to X\) is an invertible bimeasurable transformation satisfying \(\mu(T^{-1}A)=\mu(A)\) for every \(A\in\mathcal{B}\). Such a system is said to be \emph{ergodic} if every \(T\)-invariant measurable set has measure \(0\) or \(1\), and \emph{weakly mixing} if the product system \((X\times X,\mathcal{B}\otimes\mathcal{B},\mu\otimes\mu,T\times T)\) is ergodic. Equivalently, the Koopman operator \(U_T\colon L^2(X,\mu;\mathbb{C})\to L^2(X,\mu;\mathbb{C})\) defined by \(U_T f=f\circ T\) admits no nonzero complex eigenfunctions in the zero-mean subspace
\[
L_0^2(X,\mu;\mathbb{C})=\Bigl\{f\in L^2(X,\mu;\mathbb{C}):\int_X f\,d\mu=0\Bigr\}.
\]
We simply write \(L^2(\mu)\) and \(L_0^2(\mu)\) whenever the scalar
field is clear from the context.

For a set \(E\subseteq\mathbb N\), its upper Banach density is
defined by
\[
d^*(E)
:=
\limsup_{L\to\infty}\sup_{N\geq1}
\frac{|E\cap[N,N+L-1]|}{L}.
\]
A sequence \((a_n)_{n\geq1}\) is said to converge to \(a\) in
\emph{uniform density}, written
\[
\operatorname{UD}\!-\!\lim_{n\to\infty}a_n=a,
\]
if, for every \(\eta>0\), the set
\[
\{n\in\mathbb N:|a_n-a|>\eta\}
\]
has upper Banach density zero.

\begin{proposition}[Uniform-density weak mixing of all orders]
\label{prop:uniform-density-wm}
Let \((X,\mathcal B,\mu,T)\) be an invertible weakly mixing
probability-preserving system. Let \(k\geq1\) and let
\(A_0,A_1,\ldots,A_k\in\mathcal B\). Then
\[
\operatorname{UD}\!-\!\lim_{n\to\infty}
\mu\left(
A_0\cap T^{-n}A_1\cap\cdots\cap T^{-kn}A_k
\right)
=
\prod_{j=0}^k\mu(A_j).
\]
\end{proposition}

\begin{proof}
This follows from the uniform weakly mixing PET theorem (Theorem~3.1
of \cite{Ber87}), applied to the essentially distinct nonconstant
linear polynomials
\[
p_j(n)=jn,\qquad 1\leq j\leq k.
\]
Indeed, the weakly mixing PET theorem gives
\[
\operatorname{UD}\!-\!\lim_{n\to\infty}
\int_X
1_{A_0}\prod_{j=1}^k U_T^{jn}1_{A_j}\,d\mu
=
\prod_{j=0}^k\mu(A_j),
\]
where \(U_Tf=f\circ T\). Since
\[
\int_X
1_{A_0}\prod_{j=1}^k U_T^{jn}1_{A_j}\,d\mu
=
\mu\left(
A_0\cap T^{-n}A_1\cap\cdots\cap T^{-kn}A_k
\right),
\]
the result follows from \cite{Ber87}.
\end{proof}

\begin{proposition}[Uniform-density convergence and uniform Ces\`aro averages]
\label{prop:ud-to-uniform-cesaro}
Let \((a_n)_{n\geq1}\) be a bounded sequence of complex numbers. Then
the following statements are equivalent:
\begin{enumerate}
\item
\[
\operatorname{UD}\!-\!\lim_{n\to\infty}a_n=0.
\]
\item
\[
\lim_{L\to\infty}\sup_{N\geq1}
\frac1L\sum_{n=N}^{N+L-1}|a_n|=0.
\]
\end{enumerate}
\end{proposition}

\begin{proof}
Suppose first that \(a_n\to0\) in uniform density, and set
\[
M:=\sup_{n\geq1}|a_n|<\infty.
\]
For \(\eta>0\), put
\[
E_\eta:=\{n\in\mathbb N:|a_n|>\eta\}.
\]
For every \(N,L\geq1\),
\[
\frac1L\sum_{n=N}^{N+L-1}|a_n|
\leq
\eta+
M\frac{|E_\eta\cap[N,N+L-1]|}{L}.
\]
Consequently,
\[
\limsup_{L\to\infty}\sup_{N\geq1}
\frac1L\sum_{n=N}^{N+L-1}|a_n|
\leq
\eta+Md^*(E_\eta)
=\eta.
\]
Since \(\eta>0\) is arbitrary, the required limit is zero.

Conversely, for every \(\eta>0\),
\[
\eta\,
\frac{|E_\eta\cap[N,N+L-1]|}{L}
\leq
\frac1L\sum_{n=N}^{N+L-1}|a_n|.
\]
Taking the supremum over \(N\), followed by the limit superior as
\(L\to\infty\), gives \(d^*(E_\eta)=0\). Hence \(a_n\to0\) in
uniform density.
\end{proof}

Applying Propositions~\ref{prop:uniform-density-wm} and
\ref{prop:ud-to-uniform-cesaro}, we obtain, for every \(k\geq1\)
and every \(A_0,\ldots,A_k\in\mathcal B\),
\begin{equation}
\label{eq:bergelson-pet}
\lim_{L\to\infty}\sup_{N\geq1}
\frac1L\sum_{n=N}^{N+L-1}
\left|
\mu\left(
A_0\cap T^{-n}A_1\cap\cdots\cap T^{-kn}A_k
\right)
-
\prod_{j=0}^k\mu(A_j)
\right|
=0.
\end{equation}

\subsection{Spectral measures and weak mixing}

Let \((X,\mathcal B,\mu,T)\) be a probability-preserving system.
Equip \(L^2(X,\mu;\mathbb C)\) with the inner product
\[
\langle f,g\rangle
:=
\int_X f\,\overline{g}\,d\mu.
\]
Recall that the Koopman operator is
\[
U_Tf:=f\circ T.
\]
For every \(f\in L_0^2(X,\mu;\mathbb C)\), the sequence
\[
n\longmapsto\langle f,U_T^nf\rangle
\]
is positive definite. Consequently, by the Herglotz--Bochner theorem
(see \cite{CFS82,Kat04}), there exists a unique finite positive Borel
measure \(\sigma_f\) on \(\mathbb T\), called the \emph{spectral
measure} of \(f\), such that
\[
\langle f,U_T^nf\rangle
=
\int_{\mathbb T}e(-nt)\,d\sigma_f(t)
\qquad(n\in\mathbb Z).
\]

When \(f\) is real-valued, define
\[
r_f(n)
:=
\int_X f(f\circ T^n)\,d\mu
=
\langle f,U_T^nf\rangle.
\]
Then \(\sigma_f\) is symmetric, meaning that
\[
\sigma_f(-E)=\sigma_f(E)
\]
for every Borel set \(E\subseteq\mathbb T\), and
\[
\begin{aligned}
r_f(n)
&=
\int_{\mathbb T}e(-nt)\,d\sigma_f(t)\\
&=
\int_{\mathbb T}e(nt)\,d\sigma_f(t)\\
&=
\int_{\mathbb T}\cos(2\pi nt)\,d\sigma_f(t)
\in\mathbb R.
\end{aligned}
\]
In particular,
\[
r_f(-n)=r_f(n),
\qquad
r_f(0)=\sigma_f(\mathbb T)=\|f\|_2^2.
\]

A finite Borel measure on \(\mathbb T\) is called
\emph{continuous} if it is nonatomic. A probability-preserving
system is weakly mixing if and only if the spectral measure of every
function in \(L_0^2(\mu)\) is continuous \cite{CFS82}.

For a finite Borel measure \(\sigma\) on \(\mathbb T\), define
\[
\widehat{\sigma}(n)
:=
\int_{\mathbb T}e(-nt)\,d\sigma(t)
\qquad(n\in\mathbb Z).
\]
If \(\sigma\) is a symmetric probability measure, then
\[
\widehat{\sigma}(n)
=
\int_{\mathbb T}\cos(2\pi nt)\,d\sigma(t)
\in[-1,1],
\]
and
\[
\widehat{\sigma}(-n)=\widehat{\sigma}(n),
\qquad
\widehat{\sigma}(0)=1.
\]

\subsection{Stationary Gaussian systems and Hermite polynomials}
Let \(\sigma\) be a symmetric probability measure on \(\mathbb{T}\). Since \((\hat{\sigma}(n))_{n\in\mathbb{Z}}\) is a real positive-definite sequence satisfying \(\hat{\sigma}(0)=1\), one can construct a centered, real-valued stationary Gaussian process \((Z_n)_{n\in\mathbb{Z}}\) on the product probability space \((X,\mathcal{B},\mu)=(\mathbb{R}^{\mathbb{Z}},\mathcal{B}(\mathbb{R})^{\otimes\mathbb{Z}},\mu_\sigma)\) equipped with the shift transformation \(T((z_k)_{k\in\mathbb{Z}})=(z_{k+1})_{k\in\mathbb{Z}}\), coordinate projections \(Z_n((z_k)_{k\in\mathbb{Z}})=z_n=(Z_0\circ T^n)((z_k)_{k\in\mathbb{Z}})\), and covariance structure given by
\[
\mathbb{E}[Z_j Z_k]=\hat{\sigma}(j-k) \qquad (j,k\in\mathbb{Z}).
\]
The quadruple \((X,\mathcal{B},\mu,T)\) is referred to as the \emph{Gaussian system} associated with \(\sigma\) \cite[Chapter~14]{CFS82}. By the Maruyama--Fomin theorem, the Gaussian shift associated with
\(\sigma\) is ergodic if and only if \(\sigma\) is nonatomic. For
Gaussian shifts, ergodicity is equivalent to weak mixing.

Let \(\gamma\) denote the standard normal distribution \(N(0,1)\) on \(\mathbb{R}\), and let \((H_m)_{m\ge 0}\) denote the orthonormal Hermite polynomials in \(L^2(\mathbb{R},\gamma)\), given explicitly in low degrees by
\[
H_0(z)=1, \qquad H_1(z)=z, \qquad H_2(z)=\frac{z^2-1}{\sqrt{2}}, \qquad H_3(z)=\frac{z^3-3z}{\sqrt{6}}.
\]
For any pair of jointly Gaussian random variables \((U,V)\) with \(U,V\sim N(0,1)\) and \(\mathbb{E}[UV]=\rho\in[-1,1]\), Mehler's formula provides the orthogonality relation
\begin{equation}\label{eq:mehler}
\mathbb{E}\bigl[H_k(U)H_m(V)\bigr]=\delta_{km}\rho^k \qquad (k,m\ge 0).
\end{equation}

\section{Two-sided \(2\)-nice recurrence: definitions, spectral characterization, and proof}
\label{sec:characterization-proof}
Let \((X,\mathcal{B},\mu,T)\) be an invertible probability-preserving system. For a real-valued function \(g\in L^\infty(X,\mu)\) and an integer \(n\in\mathbb{Z}\), we define the three-point correlation functional
For a measurable function \(h\) and \(m\in\mathbb Z\), we use the
notation
\[
T^mh:=h\circ T^m.
\]

\[
\Lambda_n(g):=\int_X g(T^n g)(T^{2n}g)\,d\mu,
\]
and for a measurable set \(A\in\mathcal{B}\), we write
\[
b_n(A):=\mu(A\cap T^{-n}A\cap T^{-2n}A)=\Lambda_n(1_A).
\]

\begin{definition}[One-sided 2-nice recurrence]\label{def:one-sided}
A subset \(R\subseteq\mathbb{N}\) is called a \emph{set of one-sided 2-nice recurrence for weakly mixing systems} if for every weakly mixing system \((X,\mathcal{B},\mu,T)\), every \(A\in\mathcal{B}\) with \(\mu(A)>0\), and every \(\epsilon>0\), there exists \(n\in R\) such that
\[
b_n(A)\ge\mu(A)^3-\epsilon.
\]
\end{definition}

\begin{definition}[Two-sided 2-nice recurrence]\label{def:two-sided}
A subset \(R\subseteq\mathbb{N}\) is called a \emph{set of two-sided 2-nice recurrence for weakly mixing systems} if for every weakly mixing system \((X,\mathcal{B},\mu,T)\), every \(A\in\mathcal{B}\), and every \(\epsilon>0\), there exists \(n\in R\) such that
\[
b_n(A)+b_n(A^c)\ge\mu(A)^3+\mu(A^c)^3-\epsilon.
\]
Equivalently, \(R\) \emph{fails two-sided 2-nice recurrence} (or is said to be \emph{bad in the two-sided sense}) if there exist a weakly mixing system \((X,\mathcal{B},\mu,T)\), a measurable set \(A\in\mathcal{B}\), and a constant \(\epsilon_0>0\) such that
\[
b_n(A)+b_n(A^c)\le\mu(A)^3+\mu(A^c)^3-\epsilon_0 \qquad (n\in R).
\]
\end{definition}

\begin{remark}
Here ``two-sided'' refers to complement symmetrization, namely the
simultaneous treatment of \(A\) and \(A^c\); it does not refer to
positive and negative iterates.
\end{remark}

\begin{definition}[Negative spectral condition]\label{def:spectral}
We say that a subset \(R\subseteq\mathbb{N}\) satisfies the \emph{negative spectral condition} if there exist a continuous symmetric probability measure \(\sigma\) on \(\mathbb{T}\) and a constant \(\epsilon_2>0\) such that
\begin{equation}\label{eq:neg-spectral}
2\hat{\sigma}(n)+\hat{\sigma}(2n)\le-\epsilon_2 \qquad (n\in R).
\end{equation}
Dually, \(R\subseteq\mathbb{N}\) is called a \emph{set of spectral 2-nice recurrence} if for every continuous symmetric probability measure \(\sigma\) on \(\mathbb{T}\) and every \(\epsilon>0\), there exists \(n\in R\) such that
\[
2\hat{\sigma}(n)+\hat{\sigma}(2n)>-\epsilon.
\]
\end{definition}

The theorem stated below constitutes the principal result of this section. It captures a single obstruction to two-sided \(2\)-nice recurrence across three distinct formulations: the language of multiple recurrence, the language of correlation inequalities for bounded functions, and the language of Fourier coefficients of measures. This equivalence is particularly illuminating because each formulation highlights a complementary facet of the phenomenon. The dynamical statement makes explicit the failure of recurrence; the functional formulation isolates the governing quadratic correlation expression; and the spectral formulation expresses this obstruction purely in terms of a continuous symmetric measure on the circle.

Such a direct passage from a dynamical condition to a spectral criterion is not available for an arbitrary three-point correlation of a single set. Indeed, upon centering its indicator function, a three-point correlation acquires a cubic term that is not encoded by the spectral measure of the function. The two-sided formulation is purposefully designed to circumvent this obstacle: adding the corresponding correlation of the complement cancels the cubic contribution, leaving behind the specific quadratic combination \(2r_n+r_{2n}\), where \(r_k\) denotes the correlation at lag \(k\). This exact cancellation is the fundamental reason a purely spectral characterization is attainable here.

The functional condition acts as the bridge connecting these two perspectives. It records the quadratic obstruction without reference to a specific set, while the Gaussian construction together with a suitable bounded transformation enables a continuous spectral obstruction to be lifted back to a weakly mixing dynamical counterexample. Thus, the theorem is far more than a mere reformulation: it establishes a complete equivalence between the recurrence problem and a uniform negative spectral inequality. The concluding assertion further notes that this obstruction is sufficiently strong to compel the failure of ordinary one-sided \(2\)-nice recurrence.

\begin{theorem}[Three equivalent characterizations of two-sided failure]
\label{thm:main}
Let \(R\subseteq\mathbb N\). Then the following three statements are equivalent:
\begin{enumerate}
\item \textup{(Dynamical condition)} \(R\) is bad in the two-sided sense: there exist a weakly mixing system \((X,\mathcal B,\mu,T)\), a measurable set \(A\in\mathcal B\), and a constant \(\varepsilon_0>0\) such that
\[
b_n(A)+b_n(A^c)\leq \mu(A)^3+\mu(A^c)^3-\varepsilon_0
\qquad(n\in R).
\]
\item \textup{(Functional condition)} There exist an invertible weakly mixing probability-preserving system \((X,\mathcal B,\mu,T)\), a bounded real-valued function \(f\in L^\infty(X,\mu)\) with \(\int_X f\,d\mu=0\), and a constant \(\delta>0\) such that
\begin{equation}\label{eq:functional-negative-condition}
2\int_X f(T^n f)\,d\mu+\int_X f(T^{2n}f)\,d\mu
\leq-\delta\qquad(n\in R).
\end{equation}
\item \textup{(Spectral condition)} There exist a continuous symmetric probability measure \(\sigma\) on \(\mathbb T\) and a constant \(\varepsilon_2>0\) such that
\[
2\widehat\sigma(n)+\widehat\sigma(2n)\leq-\varepsilon_2
\qquad(n\in R).
\]
\end{enumerate}
Moreover, condition \textup{(3)} implies that \(R\) also fails ordinary one-sided \(2\)-nice recurrence for weakly mixing systems.
\end{theorem}

\subsection{Proof of the main spectral characterization}

Theorem~\ref{thm:main} furnishes three equivalent descriptions of the failure of two-sided \(2\)-nice recurrence: a dynamical counterexample, a bounded mean-zero functional witness, and a uniform negative inequality for the Fourier coefficients of a continuous symmetric measure. We now proceed to establish this equivalence. The proof is structured around a transparent algebraic observation: although a three-point correlation is not, in general, determined by the spectral measure of a single function, the sum of the correlations of a set and its complement is completely governed by its two-point correlations.

It is precisely at this juncture that the two-sided formulation becomes indispensable. For a single set, the centered expansion involves a cubic term that remains invisible to the ordinary spectral measure. Replacing the set by its complement reverses the sign of this cubic term while leaving the quadratic correlations invariant. Consequently, taking their sum eliminates the very component that obstructs a purely spectral description.

The proof proceeds in two directions. Starting from a dynamical witness, we extract its normalized spectral measure. Conversely, beginning with a continuous symmetric measure satisfying the requisite negative Fourier inequality, we construct a weakly mixing Gaussian system and a bounded observable whose correlations inherit this negativity. An independent Bernoulli extension then converts the observable into an actual measurable set. This two-way passage echoes in spirit the dynamical--combinatorial correspondence emphasized by Bergelson, with an exact spectral identity serving here as the decisive intermediary.

We commence with the algebraic identities that underpin this strategy.

\begin{lemma}[Triple-correlation identities]\label{lemma:identities}
Let \(g\in L^\infty(X,\mu)\) be a real-valued function, and set \(\alpha=\int_X g\,d\mu\) and \(f=g-\alpha\). For \(n\in\mathbb{Z}\), let \(r_n=\langle f,T^n f\rangle\) and \(c_n=\int_X f(T^n f)(T^{2n}f)\,d\mu\). Then for every \(n\in\mathbb{Z}\), the following identities hold:
\begin{align}
\Lambda_n(g)-\alpha^3 &= \alpha(2r_n+r_{2n})+c_n, \label{eq:id1}\\
\Lambda_n(g)+\Lambda_n(1-g)-\alpha^3-(1-\alpha)^3 &= 2r_n+r_{2n}, \label{eq:id2}\\
2r_n+r_{2n} &= \frac{1}{2}\|f+T^n f+T^{2n}f\|_2^2-\frac{3}{2}\|f\|_2^2. \label{eq:id3}
\end{align}
\end{lemma}

\begin{proof}
Expanding the product \(g(T^n g)(T^{2n}g)=(\alpha+f)(\alpha+T^n f)(\alpha+T^{2n}f)\) and integrating over \(X\), we observe that the terms linear in \(f\) vanish since \(\int_X f\,d\mu=0\). Moreover, by the \(T\)-invariance of \(\mu\), we have \(\int_X (T^n f)(T^{2n}f)\,d\mu=\int_X T^n(f\cdot T^n f)\,d\mu=\int_X f(T^n f)\,d\mu=r_n\). This immediately yields \eqref{eq:id1}.

Next, we apply \eqref{eq:id1} to the function \(1-g\), whose mean is \(1-\alpha\) and whose centered part is \(-f\). Its pair correlations remain unaltered:
\[
\int_X(-f)\,((-f)\circ T^n)\,d\mu
=\int_X f\,(f\circ T^n)\,d\mu=r_n,
\]
and likewise the correlation at lag \(2n\) is \(r_{2n}\). On the other hand, its cubic correlation undergoes a reversal of sign:
\[
\begin{aligned}
&\int_X(-f)\,((-f)\circ T^n)\,((-f)\circ T^{2n})\,d\mu\\
&\qquad=-\int_X f\,(f\circ T^n)\,(f\circ T^{2n})\,d\mu
=-c_n.
\end{aligned}
\]
Applying \eqref{eq:id1} to \(1-g\) therefore produces
\begin{equation}
\Lambda_n(1-g)-(1-\alpha)^3
=(1-\alpha)(2r_n+r_{2n})-c_n.
\label{eq:one-minus-g}
\end{equation}

Adding \eqref{eq:id1} and \eqref{eq:one-minus-g}, we obtain
\[
\begin{aligned}
&\Lambda_n(g)+\Lambda_n(1-g)
-\alpha^3-(1-\alpha)^3\\
&\quad=\alpha(2r_n+r_{2n})+c_n
 +(1-\alpha)(2r_n+r_{2n})-c_n\\
&\quad=(\alpha+1-\alpha)(2r_n+r_{2n})\\
&\quad=2r_n+r_{2n}.
\end{aligned}
\]
Here the cubic terms cancel out, while the coefficients of the quadratic correlation sum to unity. This establishes \eqref{eq:id2}, which may equivalently be written as
\[
\boxed{
\Lambda_n(g)+\Lambda_n(1-g)
=\alpha^3+(1-\alpha)^3+2r_n+r_{2n}.
}
\]

Finally, expanding the \(L^2\)-norm gives
\[
\|f+T^n f+T^{2n}f\|_2^2=3\|f\|_2^2+2\bigl(\langle f,T^n f\rangle+\langle f,T^{2n}f\rangle+\langle T^n f,T^{2n}f\rangle\bigr)=3\|f\|_2^2+2(2r_n+r_{2n}),
\]
and rearranging terms yields \eqref{eq:id3}.
\end{proof}

\begin{remark}
The crux of Lemma~\ref{lemma:identities} is that passing from \(g\) to its complement \(1-g\) preserves the quadratic correlations while inverting the sign of the cubic correlation, causing the latter to vanish upon summation. Specifically, \eqref{eq:id1} decomposes the three-point correlation into a quadratic component and a cubic remainder; \eqref{eq:id2} confirms that this remainder disappears under complement symmetrization; and \eqref{eq:id3} provides a transparent geometric interpretation of the resulting quadratic expression in terms of the squared norm of the three-term ergodic sum \(f+T^nf+T^{2n}f\). In particular, the spectral combination \(2r_n+r_{2n}\) is not an ad hoc linear combination of Fourier coefficients, but rather the exact quadratic quantity isolated by symmetrization.
\end{remark}

We are now in a position to prove Theorem~\ref{thm:main}. The forward implication is direct, whereas the reverse direction requires the Gaussian and Hermite-analytic construction.

\begin{proof}[Proof of Theorem~\ref{thm:main}]
\textbf{\((1)\implies(3)\).} Suppose that there exist a weakly mixing system \((X,\mathcal{B},\mu,T)\), a measurable set \(A\in\mathcal{B}\), and a constant \(\epsilon_0>0\) such that
\[
b_n(A)+b_n(A^c)\le\mu(A)^3+\mu(A^c)^3-\epsilon_0 \qquad (n\in R).
\]
Set \(\alpha=\mu(A)\), \(f=1_A-\alpha\), and \(r_n=\langle f,T^n f\rangle\). By identity \eqref{eq:id2} of Lemma~\ref{lemma:identities}, we have
\[
2r_n+r_{2n}=b_n(A)+b_n(A^c)-\alpha^3-(1-\alpha)^3\le-\epsilon_0 \qquad (n\in R).
\]
Observe that \(\alpha\) cannot equal \(0\) or \(1\): indeed, if \(\alpha\in\{0,1\}\), then \(b_n(A)+b_n(A^c)=\alpha^3+(1-\alpha)^3=1\), contradicting \(-\epsilon_0<0\). Consequently, \(\|f\|_2^2=\alpha(1-\alpha)>0\).

Let \(\sigma_f\) denote the spectral measure of \(f\), and define the normalized measure \(\sigma:=\sigma_f/\|f\|_2^2\). Clearly, \(\sigma\) is a symmetric probability measure on \(\mathbb{T}\). Since \(T\) is weakly mixing and \(f\in L_0^2(\mu)\), the measure \(\sigma_f\) is atomless, and hence \(\sigma\) is continuous. Because \(\hat{\sigma}(n)=r_n/\|f\|_2^2\), dividing the preceding inequality by \(\|f\|_2^2=\alpha(1-\alpha)\) yields
\[
2\hat{\sigma}(n)+\hat{\sigma}(2n)\le-\frac{\epsilon_0}{\alpha(1-\alpha)} \qquad (n\in R).
\]
Thus, the negative spectral condition holds with \(\epsilon_2=\epsilon_0/(\alpha(1-\alpha))>0\).

\medskip
\noindent\textbf{\((3)\implies(1)\).} Suppose now that \(\sigma\) is a continuous symmetric probability measure on \(\mathbb{T}\) and \(\epsilon_2>0\) satisfies
\begin{equation}\label{eq:hyp-spec}
2\hat{\sigma}(n)+\hat{\sigma}(2n)\le-\epsilon_2 \qquad (n\in R).
\end{equation}
Let \((Z_k)_{k\in\mathbb{Z}}\) be a centered stationary Gaussian process with covariance \(\mathbb{E}[Z_j Z_k]=\hat{\sigma}(j-k)\), and let \((X,\mathcal{B},\mu,T)\) be the associated Gaussian shift system. Since \(\sigma\) is continuous, the system \((X,\mathcal{B},\mu,T)\) is weakly mixing.

For \(0<\eta\le\frac{1}{2}\), consider the bounded odd function
\[
\varphi_\eta(t):=\eta\tanh(\eta t), \qquad f_\eta:=\varphi_\eta(Z_0).
\]
Evidently, \(\|f_\eta\|_\infty\le\eta\le\frac{1}{2}\) and \(\int_X f_\eta\,d\mu=0\). For each \(n\in\mathbb{Z}\), the centered Gaussian vector \((Z_0,Z_n,Z_{2n})\) has the same joint distribution as \((-Z_0,-Z_n,-Z_{2n})\). Because \(\varphi_\eta\) is odd, the product \(\varphi_\eta(Z_0)\varphi_\eta(Z_n)\varphi_\eta(Z_{2n})\) reverses sign under the transformation \((Z_0,Z_n,Z_{2n})\mapsto(-Z_0,-Z_n,-Z_{2n})\), and hence its expectation vanishes identically:
\begin{equation}\label{eq:cubic-zero}
\int_X f_\eta(T^n f_\eta)(T^{2n}f_\eta)\,d\mu=0 \qquad (n\in\mathbb{Z}).
\end{equation}

Recall that \(\gamma\) denotes the standard Gaussian probability measure on \(\mathbb R\), defined by \(\gamma(E)=\int_E\frac{1}{\sqrt{2\pi}}e^{-t^2/2}\,dt\), and \((H_k)_{k\geq 0}\) denotes the orthonormal Hermite polynomials in \(L^2(\mathbb R,\gamma)\), so that
\[
  \int_{\mathbb R} H_k(t)H_m(t)\,d\gamma(t)=\delta_{km}
  \qquad (k,m\geq 0),
\]
with \(H_1(t)=t\). Since \(Z_0\) has distribution \(\gamma\), the mapping \(h\longmapsto h(Z_0)\) preserves the \(L^2\) inner product:
\[
  \int_X h(Z_0)q(Z_0)\,d\mu
  =
  \int_{\mathbb R}h(t)q(t)\,d\gamma(t).
\]
It follows that the sequence \((H_k(Z_0))_{k\geq 0}\) forms an orthonormal system in \(L^2(X,\mu)\).

Since \(\varphi_\eta(t)=\eta\tanh(\eta t)\) is an odd function, its Hermite expansion comprises solely odd-degree terms:
\[
  f_\eta=\varphi_\eta(Z_0)
  =\sum_{j=0}^{\infty}b_{2j+1}(\eta)H_{2j+1}(Z_0),
\]
with convergence in \(L^2(X,\mu)\). Dividing by \(\eta^2\) and subtracting \(Z_0=H_1(Z_0)\), we obtain
\[
  \eta^{-2}f_\eta-Z_0
  =
  \left(\frac{b_1(\eta)}{\eta^2}-1\right)H_1(Z_0)
  +
  \sum_{j=1}^{\infty}
  \frac{b_{2j+1}(\eta)}{\eta^2}
  H_{2j+1}(Z_0).
\]
By mutual orthogonality of the summands, Parseval's identity yields
\begin{equation}\label{eq:parseval-hermite}
  \bigl\|\eta^{-2}f_\eta-Z_0\bigr\|_{L^2(\mu)}^2
  =
  \left(\frac{b_1(\eta)}{\eta^2}-1\right)^2 + \sum_{j=1}^{\infty}
  \left(\frac{b_{2j+1}(\eta)}{\eta^2}\right)^2.
\end{equation}
Observe now that \(\eta^{-2}f_\eta=\eta^{-1}\tanh(\eta Z_0)\to Z_0\) pointwise as \(\eta\to 0\), while \(|\eta^{-1}\tanh(\eta Z_0)|\le|Z_0|\in L^2(\mu)\). By the dominated convergence theorem, \(\eta^{-2}f_\eta\to Z_0\) in \(L^2(X,\mu)\) as \(\eta\to 0\). Consequently, the left-hand side of \eqref{eq:parseval-hermite} tends to zero, which forces each nonnegative term on the right-hand side to vanish in the limit. Thus,
\[
\lim_{\eta\to 0}\frac{b_1(\eta)}{\eta^2}=1 \qquad\text{and}\qquad \sum_{j=1}^{\infty}b_{2j+1}(\eta)^2=o(\eta^4)=o\bigl(b_1(\eta)^2\bigr).
\]

Moreover,
\[
b_1(\eta)
=
\int_{\mathbb R}\eta\tanh(\eta t)\,t\,d\gamma(t)>0
\qquad(\eta>0),
\]
because the integrand is nonnegative and is positive for \(t\neq0\).

Now set
\[
r_n^{(\eta)}:=\int_X f_\eta(T^n f_\eta)\,d\mu.
\]
Applying the Mehler--Hermite formula \eqref{eq:mehler}, we deduce that
\[
r_n^{(\eta)}=\sum_{j\ge 0}b_{2j+1}(\eta)^2\widehat\sigma(n)^{2j+1},
\]
and therefore
\begin{equation}\label{eq:comb-cov}
2r_n^{(\eta)}+r_{2n}^{(\eta)}=b_1(\eta)^2\bigl(2\hat{\sigma}(n)+\hat{\sigma}(2n)\bigr)+\sum_{j\ge 1}b_{2j+1}(\eta)^2\bigl(2\hat{\sigma}(n)^{2j+1}+\hat{\sigma}(2n)^{2j+1}\bigr).
\end{equation}
Because \(|\hat{\sigma}(m)|\le 1\) for every \(m\in\mathbb{Z}\), we have \(|2\hat{\sigma}(n)^{2j+1}+\hat{\sigma}(2n)^{2j+1}|\le 3\). Combining \eqref{eq:hyp-spec} and \eqref{eq:comb-cov}, we may select \(\eta\in(0,1/2]\) sufficiently small so that \(\sum_{j\ge 1}b_{2j+1}(\eta)^2\le\frac{1}{12}\epsilon_2 b_1(\eta)^2\), which gives
\begin{equation}\label{eq:neg-cov-bound}
2r_n^{(\eta)}+r_{2n}^{(\eta)}\le-\frac{3}{4}\epsilon_2 b_1(\eta)^2 \qquad (n\in R).
\end{equation}

Define \(g:=\frac{1}{2}+f_\eta\). Since \(\|f_\eta\|_\infty\le\frac{1}{2}\), we have \(0\le g\le 1\) and \(\int_X g\,d\mu=\frac{1}{2}\). By Lemma~\ref{lemma:identities}, \eqref{eq:cubic-zero}, and \eqref{eq:neg-cov-bound}, both \(g\) and \(1-g\) satisfy
\begin{equation}\label{eq:g-bound}
\Lambda_n(g)-\frac{1}{8}=\frac{1}{2}\bigl(2r_n^{(\eta)}+r_{2n}^{(\eta)}\bigr)\le-\frac{3}{8}\epsilon_2 b_1(\eta)^2 \qquad (n\in R),
\end{equation}
and the identical inequality \(\Lambda_n(1-g)-\frac{1}{8}\le-\frac{3}{8}\epsilon_2 b_1(\eta)^2\) holds for all \(n\in R\).

In the final step, we lift the \([0,1]\)-valued function \(g\) to the indicator function of a measurable set via a Bernoulli extension. Let \(\Omega=[0,1]^{\mathbb{Z}}\) be endowed with the product Lebesgue measure \(\nu\) and the bilateral shift \(\tau\). Since \((X,\mathcal{B},\mu,T)\) is weakly mixing and \((\Omega,\mathcal{F},\nu,\tau)\) is a Bernoulli shift, their product system
\[
(Y,\mathcal{C},\rho,S):=(X\times\Omega,\mathcal{B}\otimes\mathcal{F},\mu\otimes\nu,T\times\tau)
\]
is likewise weakly mixing. Consider the measurable set
\[
A:=\{(x,\omega)\in X\times\Omega:\omega_0\le g(x)\}.
\]
By Fubini's theorem, \(\rho(A)=\int_X g(x)\,d\mu(x)=\frac{1}{2}\). Moreover, for every integer \(n\ge 1\), the coordinates \(\omega_0\), \(\omega_n\), and \(\omega_{2n}\) are distinct and hence independent Uniform\([0,1]\) random variables under \(\nu\). Applying Fubini's theorem on \(X\times\Omega\) yields
\[
b_n(A)=\int_X\nu\bigl(\{\omega\in\Omega:\omega_0\le g(x),\omega_n\le g(T^n x),\omega_{2n}\le g(T^{2n}x)\}\bigr)\,d\mu(x)=\Lambda_n(g),
\]
and analogously \(b_n(A^c)=\Lambda_n(1-g)\) for all \(n\ge 1\). Adding the estimates for \(g\) and \(1-g\) from \eqref{eq:g-bound}, we arrive at
\[
b_n(A)+b_n(A^c)\le\frac{1}{4}-\frac{3}{4}\epsilon_2 b_1(\eta)^2=\rho(A)^3+\rho(A^c)^3-\epsilon_0 \qquad (n\in R),
\]
where \(\epsilon_0:=\frac{3}{4}\epsilon_2 b_1(\eta)^2>0\). This establishes statement \((1)\), while \eqref{eq:g-bound} simultaneously shows that \(R\) fails ordinary one-sided \(2\)-nice recurrence.

\medskip
\noindent\textbf{\((1)\iff(2)\).} Suppose first that \((1)\) holds. Then there exist an invertible weakly mixing system \((X,\mathcal B,\mu,T)\), a measurable set \(A\in\mathcal B\), and \(\varepsilon_0>0\) such that
\[
b_n(A)+b_n(A^c)\leq\mu(A)^3+\mu(A^c)^3-\varepsilon_0
\qquad(n\in R).
\]
Setting \(f=1_A-\mu(A)\) yields a bounded, real-valued, mean-zero function. By the complement-cancellation identity \eqref{eq:id2},
\[
2\int_X f(T^nf)\,d\mu+
\int_X f(T^{2n}f)\,d\mu
\leq-\varepsilon_0
\qquad(n\in R),
\]
so \((2)\) holds with \(\delta=\varepsilon_0\).

Conversely, assume that \((2)\) holds. Scaling \(f\) by a sufficiently small positive constant if necessary, we may assume without loss of generality that \(\|f\|_\infty\leq\frac12\), while \eqref{eq:functional-negative-condition} remains valid with a possibly smaller positive constant \(\delta\). Define
\[
g=\frac12+f,
\]
so that \(0\leq g\leq1\) and \(\int_Xg\,d\mu=\frac12\). Let \(\Omega=[0,1]^{\mathbb Z}\) with product Lebesgue measure \(\nu\) and Bernoulli shift \(\tau\), and form the weakly mixing product system
\[
(Y,\mathcal C,\rho,S)
=(X\times\Omega,\mathcal B\otimes\mathcal F,\mu\otimes\nu,T\times\tau).
\]
Defining
\[
A=\{(x,\omega)\in X\times\Omega:\omega_0\leq g(x)\},
\]
we have \(\rho(A)=\int_Xg\,d\mu=\frac12\). Since for each \(n\in\mathbb N\) the coordinates \(\omega_0,\omega_n,\omega_{2n}\) are pairwise distinct and independent, Fubini's theorem gives
\[
b_n(A)=\int_Xg(T^ng)(T^{2n}g)\,d\mu,
\qquad
b_n(A^c)=\int_X(1-g)(1-T^ng)(1-T^{2n}g)\,d\mu.
\]
Applying the complement-cancellation identity \eqref{eq:id2} to \(g\), whose mean is \(1/2\) and whose centered component is \(f\), we obtain
\[
\begin{aligned}
b_n(A)+b_n(A^c)-\rho(A)^3-\rho(A^c)^3
&=2\int_X f(T^nf)\,d\mu+
\int_X f(T^{2n}f)\,d\mu\\
&\leq-\delta
\qquad(n\in R).
\end{aligned}
\]
Hence \(A\) witnesses the failure of two-sided \(2\)-nice recurrence, establishing \((1)\).
\end{proof}

The preceding proof establishes an exact spectral characterization for the three-term arithmetic progression \(\{0,n,2n\}\). Notably, its core ingredients are by no means restricted to the linear dilation \(n\mapsto 2n\): they depend solely on complement cancellation and the spectral representation of the surviving pair correlations. We next extend this framework to general three-point configurations.

\subsection{Extension to polynomial configurations and countable abelian groups}

Throughout this subsection, \(T^m\) denotes the \(m\)-th iterate of the invertible transformation \(T\) for any integer \(m\in\mathbb Z\), so that negative exponents are unambiguously defined.

\begin{theorem}[Polynomial configurations \(\{0,p_1(n),p_2(n)\}\)]
\label{thm:poly-config}
Let \(p_1,p_2\in\mathbb Z[n]\) be distinct polynomials with zero constant term, and suppose that for every \(n\in R\subseteq\mathbb N\),
\[
p_1(n)\neq 0,\qquad p_2(n)\neq 0,\qquad (p_2-p_1)(n)\neq 0.
\]
For a probability-preserving system \((X,\mathcal B,\mu,T)\) and \(A\in\mathcal B\), define
\[
b_n^{(p_1,p_2)}(A)
:=\mu\bigl(A\cap T^{-p_1(n)}A\cap T^{-p_2(n)}A\bigr).
\]
Then the following statements are equivalent:
\begin{enumerate}
\item There exist a weakly mixing probability-preserving system \((X,\mathcal B,\mu,T)\), a measurable set \(A\in\mathcal B\), and a constant \(\varepsilon_0>0\) such that
\[
b_n^{(p_1,p_2)}(A)+b_n^{(p_1,p_2)}(A^c)
\le
\mu(A)^3+\mu(A^c)^3-\varepsilon_0
\qquad(n\in R).
\]
\item There exist a continuous symmetric probability measure \(\sigma\) on \(\mathbb T\) and a constant \(\varepsilon_2>0\) such that
\[
\widehat\sigma(p_1(n))+\widehat\sigma(p_2(n))+\widehat\sigma((p_2-p_1)(n))
\le -\varepsilon_2
\qquad(n\in R).
\]
\end{enumerate}
\end{theorem}

\begin{proof}
We first establish \((1)\implies(2)\). Set \(\alpha:=\mu(A)\) and \(f:=1_A-\alpha\), and denote the pair correlations by \(r_m:=\int_X f\,(T^m f)\,d\mu\) for \(m\in\mathbb Z\). Expanding the product
\[
(\alpha+f)\bigl(\alpha+T^{p_1(n)}f\bigr)\bigl(\alpha+T^{p_2(n)}f\bigr)
\]
and invoking the \(T\)-invariance of \(\mu\), we deduce the generalized complement-cancellation identity
\[
b_n^{(p_1,p_2)}(A)+b_n^{(p_1,p_2)}(A^c)
-\alpha^3-(1-\alpha)^3
=
r_{p_1(n)}+r_{p_2(n)}+r_{(p_2-p_1)(n)}.
\]
Indeed, the linear terms vanish because \(\int_X f\,d\mu=0\), the cubic terms cancel upon summing with the complement, and the cross-covariance simplifies to
\[
\int_X (T^{p_1(n)}f)(T^{p_2(n)}f)\,d\mu
=
\int_X f\,(T^{p_2(n)-p_1(n)}f)\,d\mu
=
r_{(p_2-p_1)(n)}.
\]
Consequently,
\[
r_{p_1(n)}+r_{p_2(n)}+r_{(p_2-p_1)(n)}
\le -\varepsilon_0
\qquad(n\in R).
\]
Just as in the proof of Theorem~\ref{thm:main}, we have \(\alpha\notin\{0,1\}\), so \(\|f\|_2^2=\alpha(1-\alpha)>0\). Letting \(\sigma_f\) be the spectral measure of \(f\) and setting \(\sigma:=\sigma_f/\|f\|_2^2\), weak mixing guarantees that \(\sigma\) is a continuous symmetric probability measure on \(\mathbb T\). Dividing the inequality above by \(\|f\|_2^2\) yields
\[
\widehat\sigma(p_1(n))+\widehat\sigma(p_2(n))+\widehat\sigma((p_2-p_1)(n))
\le -\frac{\varepsilon_0}{\alpha(1-\alpha)}
\qquad(n\in R),
\]
which proves \((2)\).

Conversely, assume that \((2)\) holds. Let \((Z_k)_{k\in\mathbb Z}\) be a centered stationary Gaussian process with covariance \(\mathbb E[Z_jZ_k]=\widehat\sigma(j-k)\), and let \(T\) be the shift transformation on the associated Gaussian probability space. Since \(\sigma\) is continuous, this Gaussian system is weakly mixing.

For \(0<\eta\le \tfrac12\), define the bounded odd observable
\[
\varphi_\eta(t):=\eta\tanh(\eta t),\qquad
f_\eta:=\varphi_\eta(Z_0).
\]
Because the centered Gaussian triple \((Z_0,Z_{p_1(n)},Z_{p_2(n)})\) is symmetric under simultaneous sign change and \(\varphi_\eta\) is odd, its cubic expectation vanishes:
\[
\int_X f_\eta\,\bigl(T^{p_1(n)}f_\eta\bigr)\bigl(T^{p_2(n)}f_\eta\bigr)\,d\mu=0
\qquad(n\in\mathbb Z).
\]
Writing the Hermite expansion \(f_\eta=\sum_{j\ge0} b_{2j+1}(\eta)H_{2j+1}(Z_0)\), the exact same \(L^2\)-approximation argument as in Theorem~\ref{thm:main} shows that
\[
\frac{b_1(\eta)^2}{\eta^4}\to 1
\qquad\text{and}\qquad
\sum_{j\ge1} b_{2j+1}(\eta)^2
=
o\bigl(b_1(\eta)^2\bigr)
\quad\text{as }\eta\to 0.
\]
By Mehler's formula, for \(r_m^{(\eta)}:=\int_X f_\eta\,(T^m f_\eta)\,d\mu\), we obtain
\[
\begin{aligned}
r_{p_1(n)}^{(\eta)}+r_{p_2(n)}^{(\eta)}+r_{(p_2-p_1)(n)}^{(\eta)}
&=
b_1(\eta)^2
\Bigl(
\widehat\sigma(p_1(n))+\widehat\sigma(p_2(n))+\widehat\sigma((p_2-p_1)(n))
\Bigr)\\
&\qquad
+
\sum_{j\ge1} b_{2j+1}(\eta)^2
\Bigl(
\widehat\sigma(p_1(n))^{2j+1}
+\widehat\sigma(p_2(n))^{2j+1}
+\widehat\sigma((p_2-p_1)(n))^{2j+1}
\Bigr).
\end{aligned}
\]
Since \(|\widehat\sigma(m)|\le1\), the tail series is bounded in modulus by \(3\sum_{j\ge1} b_{2j+1}(\eta)^2\). Choosing \(\eta>0\) sufficiently small so that \(\sum_{j\ge1} b_{2j+1}(\eta)^2\le \frac{\varepsilon_2}{12}\,b_1(\eta)^2\) and applying \((2)\), we get
\[
r_{p_1(n)}^{(\eta)}+r_{p_2(n)}^{(\eta)}+r_{(p_2-p_1)(n)}^{(\eta)}
\le
-\frac34\,\varepsilon_2\,b_1(\eta)^2
\qquad(n\in R).
\]

Now set \(g:=\frac12+f_\eta\), so that \(0\le g\le1\) and \(\int_X g\,d\mu=\frac12\). Passing to the Bernoulli product system \((Y,\mathcal C,\rho,S):=(X\times\Omega,\mathcal B\otimes\mathcal F,\mu\otimes\nu,T\times\tau)\) and defining \(A:=\{(x,\omega)\in X\times\Omega:\omega_0\le g(x)\}\), we have \(\rho(A)=\frac12\). For every \(n\in R\), our hypothesis ensures that the indices \(0, p_1(n), p_2(n)\) are pairwise distinct, whence the Bernoulli coordinates \(\omega_0, \omega_{p_1(n)}, \omega_{p_2(n)}\) are mutually independent under \(\nu\). Fubini's theorem therefore yields
\[
b_n^{(p_1,p_2)}(A)=\int_X g\,(T^{p_1(n)}g)(T^{p_2(n)}g)\,d\mu
\quad\text{and}\quad
b_n^{(p_1,p_2)}(A^c)=\int_X (1-g)\,(T^{p_1(n)}(1-g))(T^{p_2(n)}(1-g))\,d\mu.
\]
Using the vanishing of the cubic correlation together with \(\int_X g\,d\mu=\frac12\), we find
\[
b_n^{(p_1,p_2)}(A)-\frac18
=
\frac12
\Bigl(
r_{p_1(n)}^{(\eta)}+r_{p_2(n)}^{(\eta)}+r_{(p_2-p_1)(n)}^{(\eta)}
\Bigr),
\]
with the exact same identity holding for \(A^c\). Summing the two inequalities produces
\[
b_n^{(p_1,p_2)}(A)+b_n^{(p_1,p_2)}(A^c)
\le
\frac14-\frac34\,\varepsilon_2\,b_1(\eta)^2
=
\rho(A)^3+\rho(A^c)^3-\varepsilon_0
\qquad(n\in R),
\]
where \(\varepsilon_0:=\frac34\,\varepsilon_2\,b_1(\eta)^2>0\). This completes the proof of \((1)\).
\end{proof}

\begin{remark}
\label{rem:exceptional}

The nonvanishing assumptions are used in the Bernoulli lifting step,
where the coordinates indexed by
\[
0,\qquad p_1(n),\qquad p_2(n)
\]
must be pairwise distinct. If these assumptions hold only for all
sufficiently large \(n\in R\), then
Theorem~\ref{thm:poly-config} applies to the tail
\[
R':=R\cap[N_0,\infty).
\]
Without an additional argument dealing with the finitely many
exceptional elements, this does not by itself yield the corresponding
statement for the original set \(R\).
\end{remark}

We next record the analogous characterization for actions of countable abelian groups. Recall that a probability-preserving \(G\)-action \((T_g)_{g\in G}\) is called \emph{weakly mixing} if the diagonal action on \(X\times X\) is ergodic; multiple recurrence phenomena for measure-preserving actions of countable abelian groups were investigated systematically in \cite{ABB21}. For a continuous symmetric probability measure \(\sigma\) on the Pontryagin dual group \(\widehat G\), its Fourier coefficients are defined by
\[
\widehat\sigma(g):=\int_{\widehat G}\chi(g)\,d\sigma(\chi),
\]
where symmetry means invariance under inversion:
\[
\sigma(E)=\sigma(E^{-1}),
\qquad
E^{-1}:=\{\chi^{-1}:\chi\in E\}.
\]
Here ``continuous'' means nonatomic.

\begin{theorem}[Countable abelian group actions]
\label{thm:abelian-group}
Let \(G\) be a countable abelian group and let \(\phi,\psi:G\to G\) be group homomorphisms. Suppose that for every \(g\in R\subseteq G\),
\[
\phi(g)\neq 0_G,\qquad \psi(g)\neq 0_G,\qquad (\psi-\phi)(g)\neq 0_G.
\]
For a probability-preserving \(G\)-action \((T_g)_{g\in G}\) and \(A\in\mathcal B\), define
\[
b_g^{(\phi,\psi)}(A):=\mu\bigl(A\cap T_{\phi(g)}^{-1}A\cap T_{\psi(g)}^{-1}A\bigr).
\]
Then the following statements are equivalent:
\begin{enumerate}
\item There exist a weakly mixing probability-preserving \(G\)-action \((T_g)_{g\in G}\), a measurable set \(A\in\mathcal B\), and a constant \(\varepsilon_0>0\) such that
\[
b_g^{(\phi,\psi)}(A)+b_g^{(\phi,\psi)}(A^c)
\le
\mu(A)^3+\mu(A^c)^3-\varepsilon_0
\qquad(g\in R).
\]
\item There exist a continuous symmetric probability measure \(\sigma\) on \(\widehat G\) and a constant \(\varepsilon_2>0\) such that
\[
\widehat\sigma(\phi(g))+\widehat\sigma(\psi(g))+\widehat\sigma((\psi-\phi)(g))
\le -\varepsilon_2
\qquad(g\in R).
\]
\end{enumerate}
\end{theorem}

\begin{proof}
\textbf{(1)$\Rightarrow$(2).}
Suppose that (1) holds, and put
\(\alpha=\mu(A)\) and \(f=1_A-\alpha\). For \(h\in G\), write
\[
r_h:=\int_X f\,(f\circ T_h)\,d\mu.
\]
The complement-cancellation identity gives, for every \(g\in R\),
\[
r_{\phi(g)}+r_{\psi(g)}+r_{(\psi-\phi)(g)}
\leq-\varepsilon_0.
\]
As in the proof of Theorem~\ref{thm:main}, \(0<\alpha<1\), so
\(\|f\|_2^2=\alpha(1-\alpha)>0\). Let \(\sigma_f\) be the
spectral measure of \(f\) on \(\widehat G\), and set
\(\sigma=\sigma_f/\|f\|_2^2\). The measure \(\sigma\) is a symmetric
probability measure. Weak mixing of the action implies that its
spectral measures on \(L^2_0(\mu)\) are nonatomic, so \(\sigma\) is
continuous. Dividing the preceding inequality by \(\|f\|_2^2\)
proves (2), with
\(\varepsilon_2=\varepsilon_0/\bigl(\alpha(1-\alpha)\bigr)\).

\medskip
\noindent\textbf{(2)$\Rightarrow$(1).}
Suppose that (2) holds. Since \(\widehat G\) carries a nonatomic
probability measure, it is infinite, and hence \(G\) is infinite.
Let \((Z_h)_{h\in G}\) be a centered stationary real Gaussian process
with covariance
\[
\mathbb E[Z_hZ_k]=\widehat\sigma(h-k),
\qquad h,k\in G,
\]
and let \((X,\mathcal B,\mu,(T_h)_{h\in G})\) be its shift system,
with the convention \(Z_0\circ T_h=Z_h\). The Gaussian action is
weakly mixing. Indeed, its first Gaussian chaos has maximal spectral
type \(\sigma\); the spectral types on the higher Gaussian chaoses
are dominated by convolution powers of \(\sigma\). Every such
convolution power is nonatomic, since the convolution of a nonatomic
finite measure with a finite measure on a compact group is nonatomic.
As \(G\) is abelian, this excludes finite-dimensional subrepresentations
of the Koopman representation on \(L^2_0(\mu)\).

For \(0<\eta\leq\tfrac12\), define
\[
\varphi_\eta(t)=\eta\tanh(\eta t),
\qquad f_\eta=\varphi_\eta(Z_0).
\]
This is a bounded, odd, mean-zero observable, and
\(\|f_\eta\|_\infty\leq\eta\). Write its Hermite expansion as
\[
f_\eta=\sum_{j=0}^\infty b_{2j+1}(\eta)H_{2j+1}(Z_0).
\]
As in the proof of Theorem~\ref{thm:main}, dominated convergence gives
\[
\frac{b_1(\eta)}{\eta^2}\longrightarrow1,
\qquad
\sum_{j\geq1}b_{2j+1}(\eta)^2
=o\bigl(b_1(\eta)^2\bigr)
\quad(\eta\to0),
\]
and \(b_1(\eta)>0\). Put
\[
r_h^{(\eta)}:=\int_X f_\eta(f_\eta\circ T_h)\,d\mu.
\]
Mehler's formula and symmetry of \(\sigma\) imply
\[
r_h^{(\eta)}=
\sum_{j\geq0}b_{2j+1}(\eta)^2
\widehat\sigma(h)^{2j+1}.
\]
Consequently, uniformly for \(g\in R\),
\begin{align*}
&r_{\phi(g)}^{(\eta)}+r_{\psi(g)}^{(\eta)}
+r_{(\psi-\phi)(g)}^{(\eta)}\\
&=b_1(\eta)^2\Bigl(
\widehat\sigma(\phi(g))+\widehat\sigma(\psi(g))
+\widehat\sigma((\psi-\phi)(g))\Bigr)+E_\eta(g),
\end{align*}
where
\[
|E_\eta(g)|\leq
3\sum_{j\geq1}b_{2j+1}(\eta)^2.
\]
Choose \(\eta>0\) sufficiently small that
\[
3\sum_{j\geq1}b_{2j+1}(\eta)^2
\leq\frac{\varepsilon_2}{4}b_1(\eta)^2.
\]
The spectral inequality in (2) then yields
\[
r_{\phi(g)}^{(\eta)}+r_{\psi(g)}^{(\eta)}
+r_{(\psi-\phi)(g)}^{(\eta)}
\leq-\frac{3\varepsilon_2}{4}b_1(\eta)^2
\qquad(g\in R).
\]

Set \(q=\tfrac12+f_\eta\). Then \(0\leq q\leq1\) and
\(\int q\,d\mu=\tfrac12\). Let
\(\Omega=[0,1]^G\) with product Lebesgue measure \(\nu\), and let
\((\tau_h)_{h\in G}\) be the Bernoulli shift. Since \(G\) is infinite,
this action is mixing. Thus the product action
\[
(Y,\mathcal C,\rho,(S_h)_{h\in G})
=(X\times\Omega,\mathcal B\otimes\mathcal F,
  \mu\otimes\nu,(T_h\times\tau_h)_{h\in G})
\]
is weakly mixing. Define
\[
A=\{(x,\omega)\in X\times\Omega:
      \omega_{0_G}\leq q(x)\}.
\]
Then \(\rho(A)=\tfrac12\). For each \(g\in R\), the assumptions
\(\phi(g)\neq0_G\), \(\psi(g)\neq0_G\), and
\((\psi-\phi)(g)\neq0_G\) say precisely that
\(0_G,\phi(g),\psi(g)\) are pairwise distinct. The corresponding
Bernoulli coordinates are therefore independent. Conditioning on the
base point and applying Fubini's theorem gives
\[
\begin{aligned}
 b_g^{(\phi,\psi)}(A)
 &=\int_Xq\,(q\circ T_{\phi(g)})
                    (q\circ T_{\psi(g)})\,d\mu,\\
 b_g^{(\phi,\psi)}(A^c)
 &=\int_X(1-q)\,((1-q)\circ T_{\phi(g)})
                    ((1-q)\circ T_{\psi(g)})\,d\mu.
\end{aligned}
\]
The complement-cancellation identity now shows that
\[
\begin{aligned}
&b_g^{(\phi,\psi)}(A)+b_g^{(\phi,\psi)}(A^c)
-\rho(A)^3-\rho(A^c)^3\\
&\qquad=r_{\phi(g)}^{(\eta)}+r_{\psi(g)}^{(\eta)}
+r_{(\psi-\phi)(g)}^{(\eta)}\\
&\qquad\leq-\frac{3\varepsilon_2}{4}b_1(\eta)^2
\qquad(g\in R).
\end{aligned}
\]
Thus (1) holds with
\(\varepsilon_0=\tfrac34\varepsilon_2b_1(\eta)^2>0\).
\end{proof}

\begin{remark}[Relation with one-sided \(2\)-nice recurrence]
While the functional condition \eqref{eq:functional-negative-condition} completely characterizes the failure of two-sided \(2\)-nice recurrence, it does not by itself characterize the failure of ordinary one-sided \(2\)-nice recurrence. Indeed, for \(f=1_A-\mu(A)\), the one-sided deviation contains the additional cubic term
\[
\int_X f(T^nf)(T^{2n}f)\,d\mu,
\]
as exhibited in \eqref{eq:id1}, which is eliminated only after summing the correlations of \(A\) and \(A^c\).
\end{remark}

\subsection{Unconditional characterization and separation from the ergodic class}

The spectral characterizations established thus far pertain to weakly mixing systems. In this subsection, we dispense with the weak-mixing hypothesis altogether. Remarkably, the sole modification required in the spectral criterion is that the spectral measure need no longer be continuous.

\begin{definition}[Unconditional two-sided \(2\)-nice recurrence]
A subset \(R\subseteq\mathbb N\) is said to be \emph{unconditionally two-sided \(2\)-nice recurrent} if, for every invertible probability-preserving system \((X,\mathcal B,\mu,T)\), every \(A\in\mathcal B\), and every \(\varepsilon>0\), there exists \(n\in R\) such that
\[
b_n(A)+b_n(A^c)\geq \mu(A)^3+\mu(A^c)^3-\varepsilon.
\]
\end{definition}

\begin{theorem}[Unconditional two-sided spectral characterization]
\label{thm:unconditional-characterization}
Let \(R\subseteq\mathbb N\) be nonempty. The following statements are equivalent:
\begin{enumerate}
\item \(R\) is unconditionally two-sided \(2\)-nice recurrent.
\item For every symmetric Borel probability measure \(\sigma\) on \(\mathbb T\) and every \(\varepsilon>0\), there exists \(n\in R\) such that
\[
2\widehat{\sigma}(n)+\widehat{\sigma}(2n)>-\varepsilon.
\]
\item For every symmetric Borel probability measure \(\sigma\) on \(\mathbb T\),
\[
\sup_{n\in R}
\left(2\widehat{\sigma}(n)+\widehat{\sigma}(2n)\right)\ge 0.
\]
\end{enumerate}
\end{theorem}

\begin{proof}
The equivalence of \((2)\) and \((3)\) is an immediate consequence of the definition of the supremum. We establish the equivalence between \((1)\) and \((2)\) by proving that their negations are equivalent.

Suppose first that \((1)\) fails. Then there exist an invertible probability-preserving system \((X,\mathcal B,\mu,T)\), a measurable set \(A\in\mathcal B\), and a constant \(\varepsilon_0>0\) such that
\[
b_n(A)+b_n(A^c)
\le
\mu(A)^3+\mu(A^c)^3-\varepsilon_0
\qquad(n\in R).
\]
Set \(\alpha:=\mu(A)\), \(f:=1_A-\alpha\), and \(r_m(f):=\int_X f(f\circ T^m)\,d\mu\). By the complement-cancellation identity \eqref{eq:id2}, we obtain
\[
2r_n(f)+r_{2n}(f)\le-\varepsilon_0
\qquad(n\in R).
\]
In particular, \(f\neq0\), which ensures that \(\|f\|_2^2=\alpha(1-\alpha)>0\). Letting \(\sigma_f\) denote the spectral measure of \(f\) and setting \(\sigma:=\sigma_f/\|f\|_2^2\), we obtain a symmetric Borel probability measure on \(\mathbb T\) (not necessarily continuous) satisfying \(\widehat{\sigma}(m)=r_m(f)/\|f\|_2^2\). Consequently,
\[
2\widehat{\sigma}(n)+\widehat{\sigma}(2n)
\le
-\frac{\varepsilon_0}{\alpha(1-\alpha)}
\qquad(n\in R),
\]
which shows that \((2)\) fails.

Conversely, suppose that \((2)\) fails. Then there exist a symmetric Borel probability measure \(\sigma\) on \(\mathbb T\) and a constant \(\delta>0\) such that
\[
2\widehat{\sigma}(n)+\widehat{\sigma}(2n)
\le-\delta
\qquad(n\in R).
\]
Let \((Z_k)_{k\in\mathbb Z}\) be a centered stationary real Gaussian process satisfying \(\mathbb E[Z_jZ_k]=\widehat{\sigma}(j-k)\), and let \(T\) be the bilateral shift on its associated Gaussian probability space; note that no ergodicity assumption is needed here.

For \(0<\eta\le\frac12\), define \(\varphi_\eta(t):=\eta\tanh(\eta t)\) and \(f_\eta:=\varphi_\eta(Z_0)\). Then \(f_\eta\) is bounded, odd, and mean zero, with Hermite expansion
\[
f_\eta
=
\sum_{j\ge0}b_{2j+1}(\eta)H_{2j+1}(Z_0).
\]
Since \(\eta^{-2}f_\eta=\eta^{-1}\tanh(\eta Z_0)\longrightarrow Z_0\) in \(L^2\) as \(\eta\to 0\), Parseval's identity yields
\[
\frac{b_1(\eta)}{\eta^2}\longrightarrow1
\qquad\text{and}\qquad
\sum_{j\ge1}b_{2j+1}(\eta)^2
=
o\bigl(b_1(\eta)^2\bigr).
\]
Writing \(r_m^{(\eta)}:=\int_X f_\eta(T^mf_\eta)\,d\mu\), Mehler's formula gives
\[
\begin{aligned}
2r_n^{(\eta)}+r_{2n}^{(\eta)}
&=
b_1(\eta)^2
\left(2\widehat{\sigma}(n)+\widehat{\sigma}(2n)\right)
+
\sum_{j\ge1}b_{2j+1}(\eta)^2
\left(
2\widehat{\sigma}(n)^{2j+1}
+\widehat{\sigma}(2n)^{2j+1}
\right).
\end{aligned}
\]
Since the second sum is bounded in modulus by \(3\sum_{j\ge1}b_{2j+1}(\eta)^2\), we may choose \(\eta>0\) sufficiently small so that \(3\sum_{j\ge1}b_{2j+1}(\eta)^2\le \frac{\delta}{4}b_1(\eta)^2\). It follows that
\[
2r_n^{(\eta)}+r_{2n}^{(\eta)}
\le
-\frac{3\delta}{4}b_1(\eta)^2
\qquad(n\in R).
\]

Now put \(g:=\frac12+f_\eta\), so that \(0\le g\le1\) and \(\int_Xg\,d\mu=\frac12\). Passing to the Bernoulli product system \((Y,\mathcal C,\rho,S):=(X\times\Omega,\mathcal B\otimes\mathcal F,\mu\otimes\nu,T\times\tau)\) and defining \(A:=\{(x,\omega)\in Y:\omega_0\le g(x)\}\), we have \(\rho(A)=\frac12\). By the mutual independence of \(\omega_0,\omega_n,\omega_{2n}\) and the complement-cancellation identity, we conclude that
\[
\begin{aligned}
b_n(A)+b_n(A^c)-\rho(A)^3-\rho(A^c)^3
&=
2r_n^{(\eta)}+r_{2n}^{(\eta)}\\
&\le
-\frac{3\delta}{4}b_1(\eta)^2
\qquad(n\in R).
\end{aligned}
\]
Hence \(R\) fails to be unconditionally two-sided \(2\)-nice recurrent. This completes the proof.
\end{proof}

We next deduce a positive-density consequence from the weak-mixing multiple-correlation estimate \eqref{eq:bergelson-pet}.

\begin{proposition}[Positive upper Banach density]
\label{prop:positive-density-wm}
If a subset \(R\subseteq\mathbb N\) possesses positive upper Banach density, then \(R\) is two-sided \(2\)-nice recurrent for every weakly mixing system.
\end{proposition}

\begin{proof}
Suppose, for the sake of contradiction, that the conclusion fails. Then there exist a weakly mixing system, a measurable set \(A\in\mathcal B\), and a constant \(\varepsilon_0>0\) such that
\[
\left[b_n(A)-\mu(A)^3\right]
+\left[b_n(A^c)-\mu(A^c)^3\right]
\leq-\varepsilon_0\qquad(n\in R).
\]
Consequently, for every finite interval \(I\) of positive integers,
\[
\frac1{|I|}\sum_{n\in I}
\left(
 |b_n(A)-\mu(A)^3|
 +|b_n(A^c)-\mu(A^c)^3|
\right)
\geq \varepsilon_0\frac{|R\cap I|}{|I|}.
\]
Because \(R\) has positive upper Banach density, there exists a sequence of intervals \(I_j\) with \(|I_j|\to\infty\) along which the ratios \(|R\cap I_j|/|I_j|\) remain bounded away from zero. The displayed averages are therefore bounded away from zero along these intervals. This directly contradicts \eqref{eq:bergelson-pet}?Aiapplied separately to the triples \((A,A,A)\) and \((A^c,A^c,A^c)\)?Aiwhich asserts that these averages converge to zero uniformly with respect to the initial point of the interval.
\end{proof}

\begin{example}[Ergodic versus weakly mixing recurrence]
\label{ex:ergodic-wm-separation}
Consider the arithmetic progression
\[
R=\{n\in\mathbb N:n\equiv1\pmod3\}
 =\{4,7,10,\ldots\}.
\]
Its upper Banach density is \(d^*(R)=1/3>0\), and therefore Proposition~\ref{prop:positive-density-wm} guarantees that \(R\) is two-sided \(2\)-nice recurrent for every weakly mixing system. Nevertheless, \(R\) fails to be two-sided \(2\)-nice recurrent for all ergodic systems.

To verify this latter assertion, let \(X=\mathbb Z/3\mathbb Z\) be equipped with the uniform probability measure and consider the ergodic rotation \(T(x)=x+1\pmod3\). Define
\[
h(x)=\sqrt2\cos(2\pi x/3),
\qquad
 g(x)=\frac12+c h(x),
\]
where \(0<c\leq1/(2\sqrt2)\). Then \(0\leq g\leq1\), \(\int_X g\,d\mu=1/2\), and \(\int_X h(T^nh)\,d\mu=\cos(2\pi n/3)\). Setting \(f=g-1/2=ch\), we find that for every \(n\equiv1\pmod3\),
\[
2r_n(f)+r_{2n}(f)
=c^2\left(2\cos(2\pi n/3)+\cos(4\pi n/3)\right)
=-\frac32c^2.
\]
Applying the same Bernoulli extension as in the proof of Theorem~\ref{thm:main}, we obtain a measurable set \(A\) with \(\rho(A)=1/2\) whose two-sided correlation deficit is preserved. Because the product of an ergodic system with a weakly mixing Bernoulli shift is again ergodic, the resulting product system is ergodic and satisfies
\[
b_n(A)+b_n(A^c)
=\rho(A)^3+\rho(A^c)^3-\frac32c^2
\qquad(n\in R).
\]
This confirms the failure of two-sided \(2\)-nice recurrence within the ergodic class.
\end{example}

In particular, if \(\mathcal R_{\rm all}\), \(\mathcal R_{\rm erg}\), and \(\mathcal R_{\rm wm}\) denote the collections of subsets of \(\mathbb N\) that are two-sided \(2\)-nice recurrent for, respectively, all invertible probability-preserving systems, all ergodic systems, and all weakly mixing systems, then we have the chain of inclusions
\[
\mathcal R_{\rm all}\subseteq
\mathcal R_{\rm erg}\subsetneq\mathcal R_{\rm wm}.
\]
Example~\ref{ex:ergodic-wm-separation} establishes the strictness of the second inclusion, whereas whether the first inclusion is likewise strict remains an open question.

\begin{lemma}[Two-frequency lacunary Riesz product]
\label{lem:riesz-product}
Let \((n_k)_{k\ge1}\) be a sequence of positive integers satisfying the lacunarity condition
\[
\frac{n_{k+1}}{n_k}\ge 10
\qquad(k\ge1),
\]
and let \(u,v\in\mathbb R\) be real parameters such that
\[
0<|u|+|v|\le\frac12.
\]
Then there exists a continuous symmetric probability measure \(\sigma\) on \(\mathbb T\) satisfying
\[
\widehat{\sigma}(n_k)=u
\qquad\text{and}\qquad
\widehat{\sigma}(2n_k)=v
\qquad(k\ge1).
\]
\end{lemma}

\begin{proof}
Let \(m_{\mathbb T}\) denote the normalized Haar measure on \(\mathbb T\). For each integer \(K\ge1\), define the trigonometric polynomial
\[
P_K(t)
:=
\prod_{k=1}^{K}
\left(
1+2u\cos(2\pi n_kt)+2v\cos(4\pi n_kt)
\right).
\]
Each factor is nonnegative because
\[
1+2u\cos\theta+2v\cos(2\theta)
\ge
1-2|u|-2|v|
\ge0,
\]
and hence \(P_K(t)\ge0\) for all \(t\in\mathbb T\).

We first verify the requisite dissociation property. Suppose that
\[
\sum_{j=1}^{J}d_jn_j=0,
\qquad
d_j\in\{-4,-3,\ldots,3,4\}.
\]
If not all coefficients vanish, let \(k\) be the largest index for which \(d_k\neq0\). By the lacunarity hypothesis,
\[
\sum_{j<k}n_j
\le
n_k\sum_{\ell=1}^{\infty}10^{-\ell}
=
\frac{n_k}{9}.
\]
It follows that
\[
|d_k|n_k
\le
\sum_{j<k}|d_j|n_j
\le
4\sum_{j<k}n_j
\le
\frac{4}{9}n_k,
\]
which is impossible since \(|d_k|\ge1\). Thus \(d_j=0\) for every \(j\).

Consequently, every integer admits at most one representation of the form
\begin{equation}\label{eq:unique-rep}
m=\sum_{j=1}^{K}\varepsilon_jn_j,
\qquad
\varepsilon_j\in\{-2,-1,0,1,2\}.
\end{equation}
In particular, the sole representation of \(0\) is the trivial one. Hence the constant Fourier coefficient of \(P_K\) equals \(1\), meaning that \(\int_{\mathbb T}P_K\,dm_{\mathbb T}=1\). Therefore,
\[
d\sigma_K(t):=P_K(t)\,dm_{\mathbb T}(t)
\]
defines a symmetric probability measure on \(\mathbb T\).

Expanding each individual factor in complex exponentials gives
\[
1+u\,e(n_kt)+u\,e(-n_kt)
  +v\,e(2n_kt)+v\,e(-2n_kt),
\]
where \(e(t)=e^{2\pi it}\). By the uniqueness of the representation \eqref{eq:unique-rep}, we obtain for every \(1\le k\le K\),
\begin{equation}\label{eq:sigmaK-fourier}
\widehat{\sigma_K}(n_k)=u
\qquad\text{and}\qquad
\widehat{\sigma_K}(2n_k)=v.
\end{equation}

Since the space of Borel probability measures on the compact group \(\mathbb T\) is weak-\(*\) compact, there exists a subsequence \((\sigma_{K_\ell})_{\ell\ge1}\) converging in the weak-\(*\) topology to a probability measure \(\sigma\). Because each \(\sigma_K\) is symmetric, so too is \(\sigma\). Passing to the limit in \eqref{eq:sigmaK-fourier} yields
\[
\widehat{\sigma}(n_k)=u
\qquad\text{and}\qquad
\widehat{\sigma}(2n_k)=v
\qquad(k\ge1).
\]

It remains to verify that \(\sigma\) has no atoms. Set \(\theta:=\max\{|u|,|v|\}\le\frac12\). Whenever an integer \(m\) admits a representation
\[
m=\sum_{j\ge1}\varepsilon_jn_j,
\qquad
\varepsilon_j\in\{-2,-1,0,1,2\},
\]
having precisely \(q\) nonzero coefficients, the uniqueness of this representation together with \eqref{eq:sigmaK-fourier} implies
\begin{equation}\label{eq:fourier-bound}
|\widehat{\sigma}(m)|
\le
\theta^q
\le
2^{-q}.
\end{equation}
If no such representation exists, then \(\widehat{\sigma}(m)=0\).

Fix an integer \(M\ge1\), and let \(\mathcal S_M\) denote the set of integers admitting a representation with at most \(M\) nonzero coefficients. We claim that
\begin{equation}\label{eq:counting-estimate}
|\mathcal S_M\cap[-N,N]|
=
O_M\bigl((\log N)^M\bigr).
\end{equation}
Indeed, suppose \(m=\sum_j\varepsilon_jn_j\) satisfies \(|m|\le N\), and let \(k\) be the largest index with \(\varepsilon_k\neq0\). Then
\[
\begin{aligned}
|m|
&\ge
|\varepsilon_k|n_k
-
\sum_{j<k}|\varepsilon_j|n_j
\ge
n_k-2\sum_{j<k}n_j
\ge
n_k-\frac{2}{9}n_k
=
\frac79n_k,
\end{aligned}
\]
whence
\begin{equation}\label{eq:bound-nk}
n_k\le\frac97N.
\end{equation}
Because \(n_k\ge n_1\,10^{k-1}\), the number of indices \(k\) satisfying \eqref{eq:bound-nk} is \(O(\log N)\). The number of coefficient vectors possessing at most \(M\) nonzero entries supported on these indices is therefore at most
\[
\sum_{q=0}^{M}
4^q\binom{O(\log N)}{q}
=
O_M\bigl((\log N)^M\bigr),
\]
establishing \eqref{eq:counting-estimate}.

We now estimate the quadratic Wiener averages. Combining the trivial bound \(|\widehat{\sigma}(m)|\le1\) with \eqref{eq:fourier-bound}, we have
\[
\frac{1}{2N+1}
\sum_{m=-N}^{N}|\widehat{\sigma}(m)|^2
\le
\frac{|\mathcal S_M\cap[-N,N]|}{2N+1}
+
4^{-(M+1)}.
\]
By \eqref{eq:counting-estimate}, taking the limit superior as \(N\to\infty\) gives
\[
\limsup_{N\to\infty}
\frac{1}{2N+1}
\sum_{m=-N}^{N}|\widehat{\sigma}(m)|^2
\le
4^{-(M+1)}.
\]
Since \(M\ge1\) is arbitrary, we conclude that
\[
\lim_{N\to\infty}
\frac{1}{2N+1}
\sum_{m=-N}^{N}|\widehat{\sigma}(m)|^2
=0.
\]
By Wiener's lemma,
\[
\sum_{t\in\mathbb T}\sigma(\{t\})^2
=
\lim_{N\to\infty}
\frac{1}{2N+1}
\sum_{m=-N}^{N}|\widehat{\sigma}(m)|^2
=0,
\]
showing that \(\sigma\) is atomless and hence continuous.
\end{proof}

\section{The one-sided consequence and cubic-driven failure of one-sided 2-nice recurrence}\label{sec:onesided}

Although Theorem~\ref{thm:main} addresses two-sided recurrence, our Gaussian construction immediately yields a failure of ordinary one-sided \(2\)-nice recurrence as well.

\begin{corollary}
\label{cor:negative-spectral-one-sided}
Suppose that there exist a continuous symmetric probability measure \(\sigma\) on \(\mathbb T\) and a constant \(\varepsilon_2>0\) such that
\[
2\widehat{\sigma}(n)+\widehat{\sigma}(2n)
\le -\varepsilon_2
\qquad(n\in R).
\]
Then \(R\) is not a set of one-sided \(2\)-nice recurrence for weakly mixing systems.
\end{corollary}

\begin{proof}
In the proof of Theorem~\ref{thm:main}, the Gaussian construction together with the Bernoulli lifting produces a weakly mixing system \((Y,\mathcal C,\rho,S)\) and a measurable set \(A\in\mathcal C\) of measure \(\rho(A)=1/2\) satisfying
\[
b_n(A)
\le
\frac18-\frac{3\varepsilon_2}{8}b_1(\eta)^2
=
\rho(A)^3-\frac{3\varepsilon_2}{8}b_1(\eta)^2
\qquad(n\in R).
\]
Since \(b_1(\eta)>0\) for all sufficiently small \(\eta>0\), this exhibits a uniform one-sided recurrence deficit along \(R\).
\end{proof}

By identity \eqref{eq:id1} of Lemma~\ref{lemma:identities}, for \(f=1_A-\mu(A)\), we have the decomposition
\begin{equation}\label{eq:onesided-expansion}
b_n(A)-\mu(A)^3=\mu(A)\bigl(2r_n+r_{2n}\bigr)+c_n.
\end{equation}
It is therefore natural to inquire whether a single measurable set \(A\) in a weakly mixing system can fail one-sided \(2\)-nice recurrence along an infinite set \(R\subseteq\mathbb{N}\) purely on account of a negative trilinear correlation \(c_n=\int_X f(T^n f)(T^{2n}f)\,d\mu\), even while its two-point correlations satisfy \(2r_n+r_{2n}>0\) on \(R\) (so that \(A\) satisfies two-sided \(2\)-nice recurrence along \(R\) with a strictly positive margin). We now demonstrate that this phenomenon indeed takes place, by coupling an odd Hermite mode (\(H_1\)) with an even Hermite mode (\(H_2\)).

\begin{theorem}[Cubic-driven one-sided deficit with positive spectral margin]
\label{thm:cubic-obstruction}
Let \((n_k)_{k\ge1}\) be a sequence of positive integers satisfying
\[
\frac{n_{k+1}}{n_k}\ge 10
\qquad(k\ge1),
\]
and set \(R:=\{n_k:k\ge1\}\). Then there exist a weakly mixing probability-preserving system \((Y,\mathcal C,\rho,S)\), a measurable set \(A\in\mathcal C\) with \(\rho(A)=1/2\), and constants \(\delta_1,\delta_2>0\) such that, for every \(n\in R\),
\[
2\widehat{\sigma_A}(n)+\widehat{\sigma_A}(2n)
\ge \delta_1,
\]
where \(\sigma_A\) denotes the normalized spectral measure of \(f_A:=1_A-\frac12\), and
\[
b_n(A)+b_n(A^c)
\ge
\rho(A)^3+\rho(A^c)^3+\frac{\delta_1}{4},
\]
while simultaneously
\[
b_n(A)\le \rho(A)^3-\delta_2.
\]
In particular, strict positivity of the quadratic spectral expression does not by itself guarantee one-sided \(2\)-nice recurrence.
\end{theorem}

\begin{proof}
Our strategy is to construct a bounded observable whose quadratic correlation expression is strictly positive, yet whose negative cubic correlation is sufficiently pronounced to induce a one-sided recurrence deficit.

Fix the constants
\[
a:=\frac14,\qquad
b:=\frac1{32},\qquad
c:=\frac1{1024}.
\]
For every sufficiently large integer \(M\), set
\[
s_M:=\frac{c}{M^2},\qquad
u:=\frac18,\qquad
v_M:=-\frac14+s_M.
\]
For large enough \(M\), we clearly have \(0<|u|+|v_M|<\frac12\). By Lemma~\ref{lem:riesz-product}, there exists a continuous symmetric probability measure \(\sigma_M\) on \(\mathbb T\) such that
\[
\widehat{\sigma_M}(n)=u
\qquad\text{and}\qquad
\widehat{\sigma_M}(2n)=v_M
\qquad(n\in R).
\]

Let \((Z_j)_{j\in\mathbb Z}\) be the centered stationary Gaussian process associated with \(\sigma_M\), having covariance \(\mathbb E[Z_jZ_\ell]=\widehat{\sigma_M}(j-\ell)\), and let \((X_M,\mathcal B_M,\mu_M,T_M)\) denote the corresponding Gaussian shift system. Because \(\sigma_M\) is continuous, this system is weakly mixing. Moreover, for every \(n\in R\), the Gaussian triple \((Z_0,Z_n,Z_{2n})\) has covariance matrix
\[
\begin{pmatrix}
1&u&v_M\\
u&1&u\\
v_M&u&1
\end{pmatrix}.
\]
In particular, its joint distribution is completely independent of the choice of \(n\in R\).

Let \(H_1(z):=z\) and \(H_2(z):=(z^2-1)/\sqrt2\) be the first two nonconstant orthonormal Hermite polynomials, and consider their bounded truncations
\[
\psi_{1,M}(z)
:=
H_1(z)1_{\{|z|\le M\}}
\qquad\text{and}\qquad
\psi_{2,M}(z)
:=
H_2(z)1_{\{|z|\le M\}}-\beta_M,
\]
where
\[
\beta_M
:=
\int_{\mathbb R}
H_2(z)1_{\{|z|\le M\}}\,d\gamma(z)
\]
and \(\gamma\) is the standard Gaussian measure. Observe that \(\psi_{1,M}\) is odd with Gaussian mean zero, whereas \(\psi_{2,M}\) is even with Gaussian mean zero. Furthermore, for every finite \(p\ge1\), we have \(\psi_{1,M}\to H_1\) and \(\psi_{2,M}\to H_2\) in \(L^p(\gamma)\) as \(M\to\infty\); more precisely, super-exponential Gaussian tail bounds yield, for every \(p<\infty\) and every \(N>0\),
\begin{equation}\label{eq:section-four-2}
\|\psi_{1,M}-H_1\|_{L^p(\gamma)}
+
\|\psi_{2,M}-H_2\|_{L^p(\gamma)}
=
O_{p,N}(M^{-N}).
\end{equation}

Setting
\[
\eta_{1,M}:=\frac{a}{M}
\qquad\text{and}\qquad
\eta_{2,M}:=\frac{b}{M^2},
\]
we define the bounded mean-zero observable
\begin{equation}\label{eq:section-four-3}
f_M
:=
\eta_{1,M}\psi_{1,M}(Z_0)
+
\eta_{2,M}\psi_{2,M}(Z_0),
\end{equation}
and introduce for comparison the untruncated observable
\begin{equation}\label{eq:section-four-4}
\widetilde f_M
:=
\eta_{1,M}H_1(Z_0)
+
\eta_{2,M}H_2(Z_0).
\end{equation}

We begin by computing the quadratic and cubic correlations of \(\widetilde f_M\). For \(m\in\mathbb Z\), let \(\widetilde r_m:=\int_{X_M}\widetilde f_M(T_M^m\widetilde f_M)\,d\mu_M\). By Mehler's formula and the orthogonality of distinct Hermite chaoses,
\begin{equation}\label{eq:section-four-5}
\widetilde r_m
=
\eta_{1,M}^2\widehat{\sigma_M}(m)
+
\eta_{2,M}^2\widehat{\sigma_M}(m)^2.
\end{equation}
Hence, for every \(n\in R\),
\begin{equation}\label{eq:section-four-6}
\widetilde Q_M
:=
2\widetilde r_n+\widetilde r_{2n}
=
\eta_{1,M}^2(2u+v_M)
+
\eta_{2,M}^2(2u^2+v_M^2).
\end{equation}
Using the identities \(2u+v_M=s_M\) and \(2u^2+v_M^2=\frac{3}{32}-\frac{s_M}{2}+s_M^2\), we obtain
\begin{equation}\label{eq:section-four-7}
\widetilde Q_M
=
s_M\eta_{1,M}^2
+
\left(
\frac{3}{32}-\frac{s_M}{2}+s_M^2
\right)\eta_{2,M}^2.
\end{equation}
Multiplying by \(M^4\) and letting \(M\to\infty\) gives
\begin{equation}\label{eq:section-four-8}
M^4\widetilde Q_M
\longrightarrow
q
:=
ca^2+\frac{3}{32}b^2.
\end{equation}
With our specific choices of \(a,b,c\), this limiting value evaluates to
\begin{equation}\label{eq:section-four-9}
q
=
\frac1{16384}+\frac3{32768}
=
\frac5{32768}>0.
\end{equation}

We next evaluate the cubic correlation
\[
\widetilde C_M
:=
\int_{X_M}
\widetilde f_M
(T_M^n\widetilde f_M)
(T_M^{2n}\widetilde f_M)\,d\mu_M
\qquad (n\in R),
\]
which again is independent of \(n\in R\) because the distribution of \((Z_0,Z_n,Z_{2n})\) does not depend on \(n\in R\). Terms containing an odd total number of Gaussian factors vanish identically, while Wick's formula yields
\[
\begin{aligned}
\mathbb E[H_1(Z_0)H_1(Z_n)H_2(Z_{2n})] &= \sqrt2\,uv_M, &
\mathbb E[H_1(Z_0)H_2(Z_n)H_1(Z_{2n})] &= \sqrt2\,u^2, \\
\mathbb E[H_2(Z_0)H_1(Z_n)H_1(Z_{2n})] &= \sqrt2\,uv_M, &
\mathbb E[H_2(Z_0)H_2(Z_n)H_2(Z_{2n})] &= 2\sqrt2\,u^2v_M.
\end{aligned}
\]
Summing these contributions gives
\begin{equation}\label{eq:section-four-10}
\widetilde C_M
=
\sqrt2\,\eta_{1,M}^2\eta_{2,M}
\bigl(u^2+2uv_M\bigr)
+
2\sqrt2\,\eta_{2,M}^3u^2v_M.
\end{equation}
Substituting \(u^2+2uv_M=-\frac3{64}+\frac{s_M}{4}\) and \(u^2v_M=-\frac1{256}+\frac{s_M}{64}\), we find
\begin{equation}\label{eq:section-four-11}
M^4\widetilde C_M
\longrightarrow
\kappa
:=
-\frac{3\sqrt2}{64}a^2b.
\end{equation}
For \(a=1/4\) and \(b=1/32\), this simplifies to
\begin{equation}\label{eq:section-four-12}
\kappa
=
-\frac{3\sqrt2}{32768}.
\end{equation}
Combining \eqref{eq:section-four-9} and \eqref{eq:section-four-12}, we observe that
\begin{equation}\label{eq:section-four-13}
\frac q2+\kappa
=
\frac{5-6\sqrt2}{65536}<0.
\end{equation}
Thus, although the limiting quadratic expression is strictly positive, the combination \(\frac q2+\kappa\) is strictly negative.

We now pass from \(\widetilde f_M\) to its bounded truncation \(f_M\). Define
\[
r_m^{(M)}
:=
\int_{X_M}f_M(T_M^mf_M)\,d\mu_M
\qquad\text{and}\qquad
C_M
:=
\int_{X_M}
f_M(T_M^nf_M)(T_M^{2n}f_M)\,d\mu_M
\quad(n\in R),
\]
noting once more that \(C_M\) does not depend on \(n\in R\). Combining \eqref{eq:section-four-2} with the Cauchy--Schwarz inequality for pair correlations and H??lder's inequality for triple correlations, we deduce that
\begin{equation}\label{eq:section-four-14}
M^4
\left|
\bigl(2r_n^{(M)}+r_{2n}^{(M)}\bigr)-\widetilde Q_M
\right|
\longrightarrow0
\end{equation}
and
\begin{equation}\label{eq:section-four-15}
M^4|C_M-\widetilde C_M|
\longrightarrow0
\end{equation}
as \(M\to\infty\), uniformly for \(n\in R\) (since the marginal distributions are standard Gaussian and the covariance matrix of \((Z_0,Z_n,Z_{2n})\) is invariant across \(n\in R\)).

Setting \(Q_M:=2r_n^{(M)}+r_{2n}^{(M)}\) for \(n\in R\), equations \eqref{eq:section-four-8}, \eqref{eq:section-four-11}, \eqref{eq:section-four-14}, and \eqref{eq:section-four-15} imply
\begin{equation}\label{eq:section-four-16}
M^4Q_M\longrightarrow q>0
\qquad\text{and}\qquad
M^4C_M\longrightarrow\kappa.
\end{equation}
In view of \eqref{eq:section-four-13}, for all sufficiently large \(M\) we have both
\begin{equation}\label{eq:section-four-18}
Q_M>0
\qquad\text{and}\qquad
\frac12Q_M+C_M<0.
\end{equation}

Furthermore, because \(\|\psi_{1,M}\|_\infty\le M\), \(\|\psi_{2,M}\|_\infty\le \frac{M^2+1}{\sqrt2}+|\beta_M|\), and \(\beta_M\to0\), we have
\[
\|f_M\|_\infty
\le
\frac{a}{M}\|\psi_{1,M}\|_\infty
+
\frac{b}{M^2}\|\psi_{2,M}\|_\infty
\le
a+\frac{b}{\sqrt2}+o(1).
\]
Since \(a+\frac{b}{\sqrt2}=\frac14+\frac1{32\sqrt2}<\frac12\), we may fix an integer \(M\) sufficiently large so that \eqref{eq:section-four-18} holds alongside
\begin{equation}\label{eq:section-four-20}
\|f_M\|_\infty\le\frac12.
\end{equation}
Suppressing \(M\) from the notation henceforth, we put \(f:=f_M\) and \(g:=\frac12+f\), so that \(0\le g\le1\) and \(\int_{X_M}g\,d\mu_M=\frac12\).

Finally, we perform the Bernoulli lifting. Let \(\Omega:=[0,1]^{\mathbb Z}\) be endowed with product Lebesgue measure \(\nu\) and bilateral Bernoulli shift \(\tau\), and form the product system
\[
(Y,\mathcal C,\rho,S)
:=
(X_M\times\Omega,
 \mathcal B_M\otimes\mathcal F,
 \mu_M\otimes\nu,
 T_M\times\tau).
\]
Because the Gaussian system is weakly mixing and the Bernoulli shift is mixing, \((Y,\mathcal C,\rho,S)\) is weakly mixing. Defining
\[
A
:=
\{(x,\omega)\in X_M\times\Omega:
  \omega_0\le g(x)\},
\]
Fubini's theorem yields \(\rho(A)=\int_{X_M}g\,d\mu_M=\frac12\). Set \(f_A:=1_A-\frac12\).

For every \(m\neq0\), the independence of the Bernoulli coordinates \(\omega_0\) and \(\omega_m\) gives, upon conditioning on \(x\),
\begin{equation}\label{eq:section-four-21}
\langle f_A,S^mf_A\rangle
=
\int_{X_M}f(T_M^mf)\,d\mu_M
=
r_m^{(M)}.
\end{equation}
Similarly, for every \(n\in R\), the mutual independence of \(\omega_0,\omega_n,\omega_{2n}\) implies
\begin{equation}\label{eq:section-four-22}
\int_Y
f_A(S^nf_A)(S^{2n}f_A)\,d\rho
=
\int_{X_M}
f(T_M^nf)(T_M^{2n}f)\,d\mu_M
=
C_M.
\end{equation}

Since \(\rho(A)=1/2\), the centered triple-correlation identity \eqref{eq:id1} together with \eqref{eq:section-four-21} and \eqref{eq:section-four-22} yields
\begin{equation}\label{eq:section-four-23}
b_n(A)-\frac18
=
\frac12
\left(
2\langle f_A,S^nf_A\rangle
+
\langle f_A,S^{2n}f_A\rangle
\right)
+
\int_Y f_A(S^nf_A)(S^{2n}f_A)\,d\rho
=
\frac12Q_M+C_M
\qquad(n\in R).
\end{equation}
By \eqref{eq:section-four-18}, the right-hand side is strictly negative. Setting \(\delta_2:=-\left(\frac12Q_M+C_M\right)>0\), we obtain
\begin{equation}\label{eq:section-four-24}
b_n(A)
=
\frac18-\delta_2
=
\rho(A)^3-\delta_2
\qquad(n\in R).
\end{equation}

It remains to verify the spectral and two-sided bounds. Because \(f_A=1_A-1/2\), we have \(\|f_A\|_2^2=\frac14\). Letting \(\sigma_A\) denote the normalized spectral measure of \(f_A\), equation \eqref{eq:section-four-21} gives
\begin{equation}\label{eq:section-four-25}
\widehat{\sigma_A}(m)
=
\frac{\langle f_A,S^mf_A\rangle}{\|f_A\|_2^2}
=
4r_m^{(M)}
\qquad(m\neq0).
\end{equation}
Hence, for every \(n\in R\),
\begin{equation}\label{eq:section-four-26}
2\widehat{\sigma_A}(n)+\widehat{\sigma_A}(2n)
=
4\bigl(2r_n^{(M)}+r_{2n}^{(M)}\bigr)
=
4Q_M
=:
\delta_1>0.
\end{equation}
Lastly, the complement-cancellation identity \eqref{eq:id2} together with \eqref{eq:section-four-21} yields
\[
\begin{aligned}
b_n(A)+b_n(A^c)
-\rho(A)^3-\rho(A^c)^3
&=
2\langle f_A,S^nf_A\rangle
+
\langle f_A,S^{2n}f_A\rangle
=
Q_M
=
\frac{\delta_1}{4},
\end{aligned}
\]
and therefore
\[
b_n(A)+b_n(A^c)
=
\rho(A)^3+\rho(A^c)^3+\frac{\delta_1}{4}
\qquad(n\in R).
\]
Combining this with \eqref{eq:section-four-24} and \eqref{eq:section-four-26} completes the proof.
\end{proof}

\end{document}